\input amstex
\input amsppt.sty
\magnification=\magstep1
\hsize=33truecc
\vsize=22.2truecm
\baselineskip=16truept
\NoBlackBoxes
\nologo
\pageno=1
\topmatter
\TagsOnRight

\def\N{\Bbb N}
\def\Z{\Bbb Z}

\def\l{\left}
\def\r{\right}
\def\b{\bigg}

\def\({\b(}
\def\[{\b[}
\def\){\b)}
\def\]{\b]}

\def\t{\text}
\def\f{\frac}

\def\ord{\roman{ord}}

\def\se {\subseteq}

\def\sm{\setminus}

\def\bi{\binom}
\def\eq{\equiv}

\def\ls{\leqslant}
\def\gs{\geqslant}

\def\ve{\varepsilon}
\def\da{\delta}

\def\Proof{\noindent{\it Proof}}
\def\Remark{\noindent{\it Remark}}

\def\Ack{\noindent {\bf Acknowledgment}}
\hbox{Int. J. Number Theory 15\,(2019), no.\,9, 1863--1893.}
\medskip
\title Restricted sums of four squares\endtitle
\author Zhi-Wei Sun \endauthor
\affil Department of Mathematics, Nanjing University
     \\Nanjing 210093, People's Republic of China
    \\  zwsun\@nju.edu.cn
    \\ {\tt http://math.nju.edu.cn/$\sim$zwsun}
 \endaffil
\abstract We refine Lagrange's four-square theorem in new ways by imposing some restrictions involving powers of two (including $1$).
For example, we show that each $n=1,2,3,\ldots$ can be written as $x^2+y^2+z^2+w^2$ ($x,y,z,w\in\N=\{0,1,2,\ldots\})$
with $|x+y-z|\in\{4^k:\ k\in\N\}$ (or $|2x-y|\in\{4^k:\ k\in\N\}$, or $x+y-z\in\{\pm 8^k:\ k\in\N\}\cup\{0\}\se\{t^3:\ t\in\Z\}$), and
that we can write any positive integer as
$x^2+y^2+z^2+w^2\ (x,y,z,w\in\Z)$ with $x+y+2z$ (or $x+2y+2z$) a power of four. We also prove that any $n\in\N$ can be written as $x^2+y^2+z^2+2w^2$ $(x,y,z,w\in\Z)$ with $x+y+z+w$ a square (or a cube).
In addition, we pose some open conjectures for further research; for example, we conjecture that any integer $n>1$ can be written as $a^2+b^2+3^c+5^d$ with $a,b,c,d\in\N$.
\endabstract
\thanks 2010 {\it Mathematics Subject Classification}.
Primary 11E25; Secondary 11D85, 11E20, 11P05.
\newline \indent The initial version of this paper was posted to arXiv in 2017. This work is supported by the National Natural Science Foundation of China (Grant No. 11571162) and the NSFC-RFBR Cooperation and Exchange Program (Grant No. 11811530072).
\endthanks
\endtopmatter
\document

\heading{1. Introduction}\endheading

The celebrated four-square theorem (cf. [N, pp.\,5-7]) proved by J. L. Lagrange in 1770 states that any $n\in\N=\{0,1,2,\ldots\}$
can be written as $x^2+y^2+z^2+w^2$ with $x,y,z,w\in\N$. Recently, the author [S17b] found that this can be refined in various ways
by requiring additionally that $P(x,y,z,w)$ is a square, where $P(x,y,z,w)$ is a suitable polynomial with integer coefficients.
(For example, we may take $P(x,y,z,w)=x^2y^2+y^2z^2+z^2x^2$.) Here is a challenging conjecture posed by the author.

\proclaim{1-3-5 Conjecture {\rm ([S17b, Conjecture 4.3(i)])}} Any $n\in\N$ can be written as $x^2+y^2+z^2+w^2$ with $x,y,z,w\in\N$ such that $x+3y+5z$ is a square.
\endproclaim

In this paper we aim to refine Lagrange's four-square theorem in a new direction by imposing restrictions involving power of two (including $2^0=1$).
[S17b, Theorem 1.1] asserts that for any $a\in\{1,4\}$ and $m\in\{4,5,6\}$ we can write each $n\in\N$ as $ax^m+y^2+z^2+w^2$ with $x,y,z,w\in\N$.
Actually the proof in [S17b] works for a stronger result which requires additionally that $x\in\{2^k:\ k\in\N\}\cup\{0\}$.
Similarly, by modifying the proof of [S17b, Theorem 1.2(i)] slightly we get that for any $a\in\{1,2\}$ each $n\in\N$
can be written as $x^2+y^2+z^2+w^2$ with $x,y,z,w\in\N$ and $a(x-y)\in\{4^k:\ k\in\N\}\cup\{0\}$.

Now we state our first theorem.

\proclaim{Theorem 1.1} {\rm (i)} Any $n\in\Z^+=\{1,2,3,\ldots\}$ can be written as $x^2+y^2+z^2+8^k$ with $x,y,z\in\N$ and $k\in\{0,1,2\}$. Also, for each $r\in\{0,1\}$ and integer $n>r$, we can write $n^2$ as  $x^2+y^2+z^2+4^{2k+r}$ with $k,x,y,z\in\N$.

{\rm (ii)} Any $n\in\Z^+$ can be written as $x^2+y^2+z^2+w^2$ with $x,y,z,w\in\N$ such that
$x-y=2^{\lfloor\ord_2(n)/2\rfloor}$, where $\ord_2(n)$ is the $2$-adic order of $n$.

{\rm (iii)} Each $n\in\N$ not of the form $2^{6k+3}\times7\ (k\in\N)$ can be written as $x^2+y^2+z^2+w^2$ with $x,y,z,w\in\N$ and $x-y\in\{8^k:\ k\in\N\}\cup\{0\}$.
Consequently, we can write any $n\in\N$ as $x^2+y^2+z^2+w^2$ with $x,y,z,w\in\Z$ and $x+y\in\{8^k:\ k\in\N\}\cup\{0\}$.

{\rm (iv)} Any $n\in\Z^+$ can be written as $x^2+y^2+z^2+w^2$ with $x,y,z,w\in\Z$ such that
$x+y+z+w=2^{\lfloor(\ord_2(n)+1)/2\rfloor}$.

\endproclaim
\Remark\ 1.1. For integers $y$ and $z$ of the same parity, clearly
$$y^2+z^2=2\l(\f{y+z}2\r)^2+2\l(\f{y-z}2\r)^2.$$
So, the first assertion in Theorem 1.1(i) implies that any $n\in\N$ can be written as $x^2+2y^2+2z^2+8^k$ with $x,y,z\in\N$ and $k\in\{0,1,2\}$.  If $a\in\Z^+$, $4\nmid a$, and $2^2=x^2+y^2+z^2+w^2$ for some $x,y,z,w\in\N$ with $ax\in\{4^k:\ k\in\N\}$, then  $x\in\{1,2\}$ and hence $a\in\{1,2\}$.
As $\{8^k:\ k\in\N\}\cup\{0\}\se\{t^3:\ t\in\N\}$, Theorem 1.1(iii) implies a conjecture stated in [S17b, Remark 1.2].
Theorem 1.1(iv) with $\ord_2(n)=0$ was first realized
by Euler in a letter to Goldbach dated June 9, 1750 (cf. Question 37278 in {\tt MathOverFlow}).
See [S, A281494] for the number of ways to represent $n\in\Z^+$ as $x^2+y^2+z^2+w^2$ with $x,y,z,w\in\Z$ and $|x|\ls |y|\ls|z|\ls|w|$
such that $x+y+z+w=2^{\lfloor(\ord_2(n)+1)/2\rfloor}$. For example,
$$14=0^2+1^2+(-2)^2+3^2\ \ \t{with}\ 0+1+(-2)+3=2=2^{\lfloor(\ord_2(14)+1)/2\rfloor}$$
and
$$107=(-1)^2+(-3)^2+(-4)^2+9^2\ \ \t{with}\ (-1)+(-3)+(-4)+9=1=2^{\lfloor(\ord_2(107)+1)/2\rfloor}.$$
It was proved in [SS] that any $n\in\N$ can be written as $x^2+y^2+z^2+w^2$ $(x,y,z,w\in\Z)$ with $x+y+z+w$ a square.
\medskip

The author (cf. [S17b, Conjecture 4.1]) conjectured that for each $\ve\in\{\pm1\}$ any $n\in\N$ can be written as $x^2+y^2+z^2+w^2$ $(x,y,z,w\in\N$) with $2x+\ve y$ a square.
Y.-C. Sun and the author [SS] confirmed this for $\ve=1$, but the case $\ve=-1$ remains unsolved.
The author (cf. [S17b, Conjecture 4.1]) also conjectured that any $n\in\N$ can be written as $x^2+y^2+z^2+w^2$ $(x,y,z,w\in\N$) with $x+3y$ a square. Our next theorem provides an advance in this direction.

\proclaim{Theorem 1.2}
{\rm (i)} Any $n\in\Z^+$ can be written as $x^2+y^2+z^2+w^2$ with $x,y,z,w\in\N$ and $|2x-y|\in\{4^k:\ k\in\N\}$.

{\rm (ii)} Let $a\in\{1,2\}$. Then any $n\in\N$ can be written as $x^2+y^2+z^2+w^2$ with $x,y,z,w\in\N$ and $2x-y\in\{\pm a8^k:\ k\in\N\}\cup\{0\}\se\{at^3:\ t\in\Z\}$.

{\rm (iii)} If any positive integer $n\eq 9\pmod{20}$ can be written as $5x^2+5y^2+z^2$ with $x,y,z\in\Z$ and $2\nmid z$, then
any $n\in\Z^+$ can be written as $x^2+y^2+z^2+w^2\ (x,y,z,w\in\Z)$ with $x+3y\in\{4^k:\ k\in\N\}$.
\endproclaim
\Remark\ 1.2. The author [S15, Remark 1.8] conjectured that for each $n\in\N$ we can write $20n+9$
as $5x^2+5y^2+z^2$ with $x,y,z\in\Z$ and $2\nmid z$.
\medskip

The author (cf. [S17b, Conjecture 4.3(iii)-(iv)]) conjectured that each $n\in\N$ can be written as
$x^2+y^2+z^2+w^2\ (x,y,z,w\in\N)$ with $x+y-z$ (or $x-y-z$) a square. In contrast, we have the following result.

\proclaim{Theorem 1.3} {\rm (i)} Any $n\in\Z^+$ can be written as $x^2+y^2+z^2+w^2$ with $x,y,z,w\in\N$ and $|x+y-z|\in\{4^k:\ k\in\N\}$.

{\rm (ii)}  Let $a\in\{1,2\}$. Then any $n\in\N$ can be written as $x^2+y^2+z^2+w^2$ with $x,y,z,w\in\N$ and $x+y-z\in\{\pm a8^k:\ k\in\N\}\cup\{0\}
\se\{at^3:\ t\in\Z\}$.
\endproclaim
\Remark\ 1.3. We conjecture that each $n\in\Z^+$ can be written as $x^2+y^2+z^2+w^2$ with $x,y,z,w\in\N$, $x\eq y\pmod2$ and $|x+y-z|\in\{4^k:\ k\in\N\}$ (cf. [S, A299825]) but we are unable to prove this which is stronger than Theorem 1.3(i). In contrast with Theorem 1.3(ii), we conjecture that any $n\in\N$ can be written as $x^2+y^2+z^2+w^2$ with $x+y-z$ an integer cube, where $x,y,z,w\in\N$, $x\gs y\ls z$ and $x\eq y\pmod2$ (cf. [S, A282091]).
\medskip

The author [S17b] conjectured that any $n\in\N$ can be written as
$x^2+y^2+z^2+w^2$ with $x,y,z,w\in\N$ such that $P(x,y,z)$ is a square, where $P(x,y,z)$ may be any of the polynomials
$$x+y-2z,\ 2x+y-z,\ 2x-y-z,\ x+2y-2z,\ 2x-y-2z.$$
 Here we give the following result.

\proclaim{Theorem 1.4} {\rm (i)} Let $c\in\{1,2\}$. Then each $n\in\Z^+$ can be written as $x^2+y^2+z^2+w^2$ with $x,y,z,w\in\Z$ and $x+y+2z\in\{c4^k:\ k\in\N\}$.

{\rm (ii)} Any $n\in\Z^+$ can be written as $x^2+y^2+z^2+w^2$ with $x,y,z,w\in\Z$ and $x+2y+2z\in\{4^k:\ k\in\N\}$.
Also, for each $n\in\Z^+$ we can write $n^2=x^2+y^2+z^2+w^2$ with  $x,y,z,w\in\Z$ and $x+2y+2z\in\{8^k:\ k\in\N\}$.

{\rm (iii)} If $n\in\Z^+$ does not belong to $\bigcup_{k\in\N}\{2^{4k},2^{4k+3}\}$, then we can write $n$ as $x^2+y^2+z^2+w^2$ with $x,y,z,w\in\Z$ and $x+2y+2z\in\{3\times 4^k:\ k\in\N\}$.
Consequently, any integer $n>1$ can be written as $x^2+y^2+z^2+w^2$ with $x,y,z,w\in\Z$ and $x+2y+2z\in\{3\times 2^k:\ k\in\N\}$.
\endproclaim
\Remark\ 1.4. Y.-C. Sun and the author [SS, Theorem 1.2(ii)] showed that for each $d=1,2,3$ and $m=2,3$, we can write any $n\in\N$ as $x^2+y^2+z^2+w^2$ $(x,y,z,w\in\Z$)
with $x+2y+2z\in\{dt^m:\ t\in\N\}$.

\medskip

In contrast with the 1-3-5 Conjecture, we have the following curious conjecture motivated by Theorem 1.4.

\proclaim{Conjecture 1.1} {\rm (i)} Any $n\in\Z^+$ can be written as $x^2+y^2+z^2+w^2$ with $x,y,z,w\in\N$ such that
$x+2(y-z)$ is a power of four (including $4^0=1$). Also, for each $n\in\Z^+$ we can write $n^2$ as $x^2+y^2+z^2+w^2$ with $x,y,z,w\in\N$
and $x+2(y-z)\in\{8^k:\ k\in\N\}$.

{\rm (ii)} Let $c\in\{1,2,4\}$. Then each $n\in\N$ can be written as $x^2+y^2+z^2+w^2$ with $x,y,z,w\in\N$ and $y\ls z$ such that
$c(2x+y-z)\in\{8^k:\ k\in\N\}\cup\{0\}\se\{t^3:\ t\in\N\}$.
\endproclaim
\Remark\ 1.5. (i) We have verified the first assertion in part (i) for all $n=1,\ldots,2\times10^7$, and Qing-Hu Hou extended the verification for $n$ up to $10^9$.
See [S, A279612 and A279616] for related data. For example,
$$111 = 9^2 + 1^2 + 5^2 + 2^2\ \  \t{with}\  9 + 2\times1 - 2\times 5 = 4^0. $$

(ii) We have verified part (ii) of Conjecture 1.1 for all $n=0,\ldots,2\times10^6$. See [S, A284343] for related data.
For example,
$$2976 = 20^2 + 16^2 + 48^2 + 4^2\ \  \t{with}\ 16 < 48\ \t{and}\ 2\times20 + 16 - 48 = 8. $$
\medskip

In 1917 Ramanujan [R] listed 55 possible quadruples $(a,b,c,d)$ of positive integers with $a\ls b\ls c\ls d$
such that any $n\in\N$ can be written as $ax^2+by^2+cz^2+dw^2$ with $x,y,z,w\in\Z$, and
54 of them were later confirmed by Dickson [D27] with the remaining one wrong (see also [W]).

\proclaim{Theorem 1.5}   {\rm (i)} Any $n\in\Z^+$ can be written as $x^2+y^2+z^2+2w^2$ with $x,y,z,w\in\N$ and $x-y=1$.
Also, any $n\in\Z^+$ can be written as $x^2+2y^2+2z^2+2w^2$ with $x,y,z,w\in\Z$ and $x+y+z=1$.

{\rm (ii)} Any $n\in\Z^+$ can be written as $x^2+y^2+z^2+2w^2$ $(x,y,z,w\in\Z)$ with $x+y+2z=1$. Also,
any $n\in\Z^+$ can be written as $x^2+4y^2+z^2+2w^2$ $(x,y,z,w\in\Z)$ with $x+2y+2z=1$, and
any $n\in\Z^+$ can be written as $x^2+2y^2+2z^2+2w^2$ $(x,y,z,w\in\Z)$ with $x+y+3z=1$.

{\rm (iii)} Any $n\in\Z^+$ can be written as $x^2+y^2+z^2+2w^2$ $(x,y,z,w\in\Z)$ with $y+z+2w=1$.
Also, any $n\in\Z^+$ can be written as $x^2+y^2+4z^2+2w^2$ $(x,y,z,w\in\Z)$ with $y+2z+2w=1$, and any $n\in\Z^+$
can be written as $x^2+2y^2+2z^2+2w^2$ $(x,y,z,w\in\Z)$ with $x+y+z+2w=1$.

{\rm (iv)} Any $n\in\Z^+$ can be written as $x^2+y^2+z^2+2w^2$ $(x,y,z,w\in\Z)$ with $y+z+w=1$.

{\rm (v)} Any $n\in\Z^+$ can be written as $x^2+y^2+2z^2+5w^2$ with $x,y,z,w\in\Z$ and $y+w=1$.
For each $\da=1,2$, any positive integer $n\not\eq2\pmod3$ can be written as $x^2+y^2+2z^2+(6/\da)w^2$ with $y+\da w=1$.

{\rm (vi)} Any $n\in\Z^+$ can be written as $x^2+y^2+z^2+3w^2$ $(x,y,z,w\in\Z)$ with $x+y+2z=2$.

{\rm (vii)} Any integer $n>4$ can be written as $x^2+y^2+z^2+2w^2$ with $x,y,z,w\in\Z$ such that $x+y+z=t^2$ for some $t=1,2$.

{\rm (viii)} Any integer $n>7$ can be written as $x^2+y^2+z^2+2w^2$ $(x,y,z,w\in\Z)$ with $x+y+2z=2t^2$ for some $t=1,2$.
\endproclaim
\Remark\ 1.6. We can prove several other results similar to the first assertion in Theorem 1.5(v).
As a supplement to the second assertion in Theorem 1.5(v), we conjecture that any positive integer $n\eq2\pmod3$
also can be written as $x^2+y^2+2z^2+6w^2$ with $x,y,z,w\in\Z$ and $y+w=1$, and that $11$ is the only positive integer
which cannot be  written as
$x^2+y^2+2z^2+3w^2$ with $x,y,z,w\in\Z$ and $y+2w=1$.
\medskip

\proclaim{Theorem 1.6} Let $m\in\{2,3\}$.

{\rm (i)}  Any $n\in\N$ can be written as $x^2+y^2+z^2+2w^2$ $(x,y,z,w\in\Z)$ with $x+y+z+w$ an $m$-th power.

{\rm (ii)} Any $n\in\N$ can be written as $x^2+y^2+z^2+2w^2$ $(x,y,z,w\in\Z)$ with $x+2y+2z$ an
  $m$-th power.
\endproclaim
\Remark\ 1.7. Our proof of Theorem 1.6 depends heavily on a new identity similar to Euler's four-square identity.
We even conjecture that any $n\in\N$ can be written as $x^2+y^2+z^2+2w^2$ $(x,y,z\in\N$ and $w\in\Z)$
with $x+y-z+w\in\{0,1\}$.
\medskip

We will prove Theorems 1.1-1.4 and 1.5-1.6 in Sections 2 and 3 respectively, and pose more related conjectures in Sec. 4.

\heading{2. Proofs of Theorems 1.1-1.4}\endheading

It is known that the set
$$E(a,b,c):=\{n\in\N:\ n\not= ax^2+by^2+cz^2\ \t{for all}\ x,y,z\in\Z\}\tag2.1$$
is infinite for any $a,b,c\in\Z^+$.

\proclaim{Lemma 2.1 {\rm (The Gauss-Legendre Theorem)}} We have
$$E(1,1,1)=E_0:=\{4^k(8l+7):\ k,l\in\N\}.\tag2.2$$
\endproclaim
\Remark\ 2.1. This is a well-known result on sums of three squares, see, e.g., [N, p.\,23] or [MW, p.\,42]).

\proclaim{Lemma 2.2} Let $m,n\in\N$ with $16\mid m$ and $16\nmid n$.
 Then $\{n-1,n-m\}\not\se E_0$, where $E_0$ is the set defined in $(2.2)$.
\endproclaim
\Proof. Clearly one of $n-1$ or $n-m$ is odd. If $n-1\eq 7\pmod 8$, then $n-m\eq n\eq 8\pmod{16}$.
If $n-m\eq 7\pmod 8$, then $n-1\eq n-m-1\eq6\pmod 8$. Thus $\{n-1,n-m\}\not\se E_0$ as desired. \qed
\medskip

\Remark\ 2.2. It follows from Lemmas 2.1 and 2.2 that any $n\in\Z^+$ can be written as $x^2+y^2+z^2+w^2$ with $x$ a power of two and $y,z,w\in\N$.
A stronger result given in [S17b, Theorem 1.2(v)] states that any positive integer can be written as $4^k(1+4x^2+y^2)+z^2$ with $k,x,y,z\in\N$.

\medskip
\noindent{\it Proof of Theorem 1.1}. (i) (a) If $n-1\not\in E(1,1,1)$, then $n=8^0+x^2+y^2+z^2$ for some $x,y,z\in\N$.

Below we assume that $n-1\in E(1,1,1)$. Then $n-1=4^k(8l+7)$ for some $k,l\in\N$.

If $k>0$, then $n-8\eq n \eq1\pmod 4$ and hence $n-8=x^2+y^2+z^2$ for some $x,y,z\in\N$.

Now we consider the case $k=0$. In this case $n=8l+8$. If $l$ is odd, then by Lemma 2.1 we can write $n-8=2^3l$ as $x^2+y^2+z^2$ with $x,y,z\in\N$.
Clearly,
$$8=0^2+0^2+0^2+8,\ 3\times 8=4^2+0^2+0^2+8,\ 5\times 8=4^2+4^2+0^2+8,\ 7\times8=4^2+4^2+4^2+8.$$
If $l$ is even and at least $8$, then $n-64=8(l+1)-64=2^3(l-7)\not\in E(1,1,1)$ and hence $n=8^2+x^2+y^2+z^2$ for some $x,y,z\in\N$.
This proves the first assertion in Theorem 1.1(i).

(b) Let $n=4^km$ with $k\in\N$, $m\in\Z^+$ and $4\nmid m$. If $m=1$, then
$n^2=0^2+0^2+0^2+4^{2k}$. If $m=1$ and $n>1$, then
$$n^2=(2^{2k-1})^2+(2^{2k-1})^2+(2^{2k-1})^2+4^{2(k-1)+1}.$$

Now let $m>1$. If $m<4$, then $m^2-16^\da\not\in E_0$ with $\da=0$.
If $m>4$, then by Lemma 2.2, for some $\da\in\{0,1\}$ we have
$m^2-16^{\da}\not\in E_0$. So $m^2-16^\da=x^2+y^2+z^2$ for some $x,y,z\in\N$ and hence
$$n^2=16^km^2=(4^kx)^2+(4^ky)^2+(4^kz)^2+4^{2(k+\da)}.$$
This proves the second assertion in Theorem 1.1(i) for $r=0$.

Finally, we handle the case $r=1$.
If $2\nmid m$, then $m^2-4\not\in E_0$. If $2\mid m$ and $m<8$, then $(m/2)^2-16^{\da}\not\in E_0$ with $\da=0$.
If $2\mid m$ and $m>8$, then $(m/2)^2-16^{\da}\not\in E_0$ for some $\da\in\{0,1\}$ (by Lemma 2.2).
Anyway, for some $\da\in\{0,1\}$ we can write $m^2-4\times16^\da$ as $x^2+y^2+z^2$ with $x,y,z\in\N$, and thus
$$n^2=16^km^2=(4^kx)^2+(4^ky)^2+(4^kz)^2+4^{2(k+\da)+1}$$
as desired.

(ii) Write $n=2^am$ with $a\in\N$, $m\in\Z^+$ and $2\nmid m$. Let $n_0=n/4^{\lfloor a/2\rfloor}=2^{a_0}m$, where
$a_0$ is $0$ or $1$ according as $a$ is even or odd. As $4\nmid n_0$, we have $2n_0-1\not\eq 7\pmod 8$.
By Lemma 2.1, $2n_0-1=u^2+v^2+(2y+1)^2$ for some $u,v,y\in\N$ with $u\eq v\pmod 2$. Let $z=(u+v)/2$ and $w=|u-v|/2$. Then
$2n_0-2=4y^2+4y+(z+w)^2+(z-w)^2$ and hence $n_0=x^2+y^2+z^2+w^2$ with $x=y+1$. It follows that
$$n=\l(2^{\lfloor a/2\rfloor}x\r)^2+\l(2^{\lfloor a/2\rfloor}y\r)^2+\l(2^{\lfloor a/2\rfloor}z\r)^2+\l(2^{\lfloor a/2\rfloor}w\r)^2$$
with
$$2^{\lfloor a/2\rfloor}x-2^{\lfloor a/2\rfloor}y=2^{\lfloor a/2\rfloor}.$$

(iii) Obviously, for any $k\in\N$ we have
$$2^{6k+3}\times7=(8^k\times6)^2+(8^k\times4)^2+(8^k\times2)^2+0^2\ \ \t{with}\ 8^k\times6+8^k\times2=8^{k+1}.$$
So it suffices to prove the first assertion in Theorem 1.1(iii) by induction.

For each $n\in\{0,1,\ldots,63\}\sm\{56\}$, we can verify via a computer that $n$ can be written as $x^2+y^2+z^2+w^2$
with $x,y,z,w\in\N$ and $x-y\in\{8^k:\ k\in\N\}\cup\{0\}$.

Now let $n\gs 64$ be an integer not of the form $2^{6k+3}\times7=64^k\times56\ (k\in\N)$,
and assume that any $m\in\{0,1,\ldots,n-1\}$ not of the form $2^{6k+3}\times7\ (k\in\N)$ can be
written as $x^2+y^2+z^2+w^2$
with $x,y,z,w\in\N$ and $x-y\in\{8^k:\ k\in\N\}\cup\{0\}$.

If $64\mid n$, then by the induction hypothesis we can write $n/64$ as $x^2+y^2+z^2+w^2$
with $x,y,z,w\in\N$ and $x-y\in\{8^k:\ k\in\N\}\cup\{0\}$, and hence
$$n=(8x)^2+(8y)^2+(8z)^2+(8w)^2\ \ \t{with}\ 8x-8y=8(x-y)\in\{8^k:\ k\in\N\}\cup\{0\}.$$

Below we suppose that $64\nmid n$.

{\it Case} 1. $n\not\in 4^k(16l+14)$ for any $k,l\in\N$.

In this case, we have $2n\not\in E(1,1,1)$ by (2.2), and hence
$2n=(2y)^2+z^2+w^2$ for some $y,z,w\in\N$ with $z\eq w\pmod 2$.
Thus $$n=y^2+y^2+\l(\f{z+w}2\r)^2+\l(\f{z-w}2\r)^2\ \t{with}\ y-y=0.$$

{\it Case} 2. $n=16l+14$ for some $l\in\N$.

In this case, $2n-1\eq3\pmod 8$ and hence by (2.2) we can write $2n-1$ as
$(2y+1)^2+z^2+w^2$ with $y,z,w\in\N$ and $z\eq w\pmod 2$. It follows that
 $$n=(y+1)^2+y^2+\l(\f{z+w}2\r)^2+\l(\f{z-w}2\r)^2\ \t{with}\ (y+1)-y=1.$$

{\it Case} 3. $n=4^k(16l+14)$ for some $k\in\{1,2\}$ and $l\in\N$.

In this case, $2n-64=4^{k+1}(8l+3k)$. In light of (2.2), $8l+3k=x^2+y^2+z^2$
for some integers $x\gs y\gs z\gs 0$. If $k=1$, then $8l+3k=(2n-64)/16\gs64/16=4$.
If $k=2$, then $8l+3k\gs 3k=6$. So, $x\gs 2$ and hence $2^{k+1}x=2v+8$
for some $v\in\N$. Therefore
$$2n-64=(2^{k+1}x)^2+(2^{k+1}y)^2+(2^{k+1}z)^2=(2v+8)^2+2(2^ky+2^kz)^2+2(2^ky-2^kz)^2$$
and hence
$$n=(v+8)^2+v^2+(2^k(y+z))^2+(2^k(y-z))^2\ \ \t{with}\ (v+8)-v=8^1.$$

The induction proof of Theorem 1.1(iii) is now complete.

(iv) We distinguish two cases.
\medskip

{\it Case} 1. $4\nmid n$.

 Let
$$\da_n=1-\ord_2(n)=\cases1&\t{if}\ 2\nmid n,
\\0&\t{if}\ 2\|n.
\endcases$$
Then
$$4^{\da_n}n-1\eq\cases3\pmod8&\t{if}\ 2\nmid n,\\1\pmod4&\t{if}\ 2\|n.\endcases$$
By Lemma 2.1, there are $u,v,w\in\Z$ such that
$$4^{\da_n}n-1=\l(2^{\da_n}u-1\r)^2+\l(2^{\da_n}v-1\r)^2+\l(2^{\da_n}w-1\r)^2.$$
If $2\|n$, then $\da_n=0$ and $u+v+w\eq n\eq0\pmod2$. In the case $2\nmid n$, since $(2u-1)^2=(2(1-u)-1)^2$, without loss of generality we may also
assume that $u+v+w\eq0\pmod 2$. Set
$$x=\f{u+v-w}2,\ y=\f{u-v+w}2,\ z=\f{-u+v+w}2.$$
Then $u=x+y$, $v=x+z$ and $w=y+z$. Therefore
$$\align 4^{\da_n}n-1=&\l(2^{\da_n}(x+y)-1\r)^2+\l(2^{\da_n}(x+z)-1\r)^2+\l(2^{\da_n}(y+z)-1\r)^2
\\=&\l(2^{\da_n}(x+y+z)-2\r)^2+\l(2^{\da_n}x\r)^2+\l(2^{\da_n}y\r)^2+\l(2^{\da_n}z\r)^2-1
\endalign$$
and hence
$$n=\l((x+y+z)-2^{1-\da_n}\r)^2+x^2+y^2+z^2=x^2+y^2+z^2+\l(2^{1-\da_n}-x-y-z\r)^2$$
with $x+y+z+(2^{1-\da_n}-x-y-z)=2^{\lfloor(\ord_2(n)+1)/2\rfloor}$.

\medskip
{\it Case} 2. $4\mid n$.

Write $n=4^kn_0$ with $k,n_0\in\Z^+$ and $4\nmid n_0$. By the above, there are $x_0,y_0,z_0,w_0\in\Z$ such that
$$n_0=x_0^2+y_0^2+z_0^2+w_0^2\ \ \t{with}\ \ x_0+y_0+z_0+w_0=2^{\lfloor(\ord_2(n_0)+1)/2\rfloor}.$$
It follows that
$$n=\l(2^kx_0\r)^2+\l(2^ky_0\r)^2+\l(2^kz_0\r)^2+\l(2^kw_0\r)^2$$
with
$$2^kx_0+2^ky_0+2^kz_0+2^kw_0=2^{k+\lfloor(\ord_2(n_0)+1)/2\rfloor}=2^{\lfloor(\ord_2(n)+1)/2\rfloor}.$$

Combining the above, we have proved Theorem 1.1. \qed

\proclaim{Lemma 2.3} {\rm (i) (Dickson [D39, pp.\,112-113])} We have
$$E(1,5,5)=\{n\in\N:\ n\eq2,3\pmod5\}\cup E_0.\tag2.3$$

{\rm (ii) ([S17a, Lemma 2.1])} Let $u$ and $v$ be integers with $u^2+v^2$ a positive multiple of $5$. Then
$u^2+v^2=x^2+y^2$ for some $x,y\in\Z$ with $5\nmid xy$.
\endproclaim

\proclaim{Lemma 2.4} Let $n\gs4$ be an integer not divisible by $64$.

{\rm (i)} We have $\{n,n-4\}\not\se E_0$.

{\rm (ii)} If $16\nmid n$, then $\{n,n-1\}\not\se E_0.$
If $16\mid n$, then $\{n,n-1,n-64\}\not\se E_0.$
\endproclaim
\Proof. (i) If $n=4^k(8l+7)$ for some $k\in\N$ and $l\in\N$, then $k<3$ as $64\nmid n$, and
hence $n-4\not\in E_0$ since
$$n-4=\cases 8l+7-4=8l+3&\t{if}\ k=0,
\\4(8l+7)-4=4(8l+6)&\t{if}\ k=1,
\\4^2(8l+7)-4=4(8(4l+3)+3)&\t{if}\ k=2.
\endcases$$

(ii) If $n=8l+7$ for some $l\in\N$, then $n-1=8l+6\not\in E_0$.
If $n=4(8l+7)$ for some $l\in\N$, then $n-1=32l+27=8(4l+3)+3\not\in E_0$. So $\{n,n-1\}\not\se E_0$ if $16\nmid n$.

Now we consider the case $16\mid n$. As $64\nmid n$, if $n\not\in E_0$, then $n=4^2(8l+7)$ for some $l\in\N$,
and hence $n-64=4^2(8l+3)\not\in E_0.$

The proof of Lemma 2.4 is now complete. \qed

\medskip
\noindent{\it Proof of Theorem 1.2}. (i) For $n=1,\ldots,15$ we can easily verify the desired result.

Now fix an integer $n\gs 16$ and assume that each $m=1,\ldots,n-1$ can be written as $x^2+y^2+z^2+w^2$ $(x,y,z,w\in\N$)
with $|2x-y|\in\{4^k:\ k\in\N\}$.

Let's first consider the case $16\mid n$. By the the induction hypothesis, there are $x,y,z,w\in\N$ for which $n/16=x^2+y^2+z^2+w^2$
with $|2x-y|\in\{4^k:\ k\in\N\}$, and hence $|2(4x)-4y|=4|2x-y|\in\{4^k:\ k\in\N\}$.

Now we suppose that $16\nmid n$. Then $\{5n-1,5n-16\}\not\se E_0$ by Lemma 2.2.
Let $\da=0$ if $5n-1\not\in E_0$, and $\da=1$ otherwise. In view of Lemma 2.1,
we can write $5n-16^\da$ as $x^2+y^2+z^2$ with $x,y,z\in\Z$.
Since a square is congruent to one of $0,1,-1$ modulo $5$.
one of $x^2,y^2,z^2$ must be congruent to $-1$ modulo $5$.
Without loss of generality, we may assume that $x^2+1\eq y^2+z^2\eq0\pmod 5$.
If $y^2+z^2\not=0$, then by Lemma 2.3(ii) we can write $y^2+z^2=y_1^2+z_1^2$ with $y_1,z_1\in\Z$ and $5\nmid y_1z_1$.
Without loss of generality, we simply assume that $x\eq  -2\times 4^\da\pmod 5$
(otherwise we use $-x$ instead of $x$)
and that either $y=z=0$ or $y\eq 2z\eq -2\times 4^{\da}\pmod5$.

Clearly, $r=(x+2\times4^\da )/5,\ s=(2x-4^\da)/5$, $u=(2y+z)/5$ and $v=(2z-y)/5$ are all integers.
Observe that
$$r^2+s^2+u^2+v^2=\f{(4^\da)^2+x^2}5+\f{y^2+z^2}5=n \ \ \t{with}\ \  2r-s=4^\da.$$
If $x>-2$, then $x\gs2$ since $x\eq(-1)^{\da-1}2\pmod5$, and hence $r,s\in\N$. When $x\ls-2$, clearly $s<0$, and
$$r>0\iff(\da=1\ \t{and}\ x=-3).$$
If $r\ls0$ and $s\ls 0$, then
$$|2|r|-|s||=|2(-r)-(-s)|=|2r-s|=4^\da.$$

Now it remains to consider the case
$\da=1$ and $x=-3$. Note that
$$\f{y^2+z^2}5=u^2+v^2=n-r^2-s^2=n-1^2-(-2)^2>0$$
and hence $y\eq 2z\eq-2\times 4^{\da}\pmod5$.
When $y\not=-3$, in the spirit of the above arguments, all the four numbers
$$\bar r=\f{y+2\times 4^\da}5,\ \bar s=\f{2y-4^\da}5,\ \bar u=\f{2x+z}5,\ \bar v=\f{2z-x}5$$
are integral, and
$$\bar r^2+\bar s^2+\bar u^2+\bar v^2=\f{(4^{\da})^2+y^2}5+\f{x^2+z^2}5=n$$
with $|2|\bar r|-|\bar s||=|2\bar r-\bar s|=4^\da$.
If $y=-3$, then
$$5n-16=x^2+y^2+z^2=(-3)^2+(-3)^2+z^2\eq z^2+2\not\eq0,1\pmod4,$$
thus $5n-1\not\eq0,3\pmod4$ and hence $5n-1\not\in E_0$ which contradicts $\da=1$.
This concludes our induction proof of Theorem 1.2(i).
\medskip

(ii) For every $n=0,1,\ldots,63$, we can verify the desired result directly.

Now fix an integer $n\gs64$ and assume that the desired result holds for all smaller values of $n$.

If $64\mid n$, then by the induction hypothesis we can write $n/64$ as $x^2+y^2+z^2+w^2$
with $x,y,z,w\in\N$ and $2x-y\in\{\pm a8^k:\ k\in\N\}\cup\{0\}$, and hence
$$n=(8x)^2+(8y)^2+(8z)^2+(8w)^2\ \ \t{with}\ 2(8x)-(8y)\in\{\pm a8^k:\ k\in\N\}\cup\{0\}.$$

Below we suppose that $64\nmid n$.

By Lemma 2.4, $5n-(a\da)^2\not\in E_0$ for some $\da\in\{0,1,8\}$ satisfying
$$\da=8\ \Longrightarrow\ (a=1\ \t{and}\ 16\mid n).\tag2.4$$
In view of (2.2),
$5n-(a\da)^2=x^2+y^2+z^2$ for some $x,y,z\in\Z$. Since any square is congruent to one of $0,\pm1$ modulo $5$,
one of $x^2,y^2,z^2$, say $x^2$, is congruent to $-(a\da)^2$ modulo $5$. Without loss of generality,
we simply suppose that $x\eq-2a\da\pmod 5$. As $y^2\eq(2z)^2\pmod 5$, we may also suppose that $y\eq2z\pmod 5$.
Thus all the numbers
$$r=\f{x+2a\da}5,\ s=\f{2x-a\da}5,\ u=\f{2y+z}5,\ v=\f{y-2z}5$$
are integral. Note that
$$n=\f{(a\da)^2+x^2}5+\f{y^2+z^2}5=r^2+s^2+u^2+v^2$$
with $$2r-s=a\da\in\{a8^k:\ k\in\N\}\cup\{0\}.$$
If $x\gs a\da/2$, then $r\gs0$ and $s\gs0$.
If $x\ls-2a\da$, then $r\ls 0$ and $s\ls-a\da\ls0$.

Now we handle the remaining case $-2a\da<x<a\da/2$.
Clearly, $r>0>s$, $a\da=2r-s>2$, hence $\da=8$, and $a=1$ and $16\mid n$ by (2.4).
As $2r-s=8$, we must have
$$(r,s)\in\{(1,-6),\ (2,-4),\ (3,-2)\}.$$
Note that $2\times2-4=0$ and $2\times2-3=1$.
If $(r,s)=(1,-6)$, then
$$n=1^2+(-6)^2+u^2+v^2\not\eq0\pmod4,$$
which contradicts $16\mid n$.
This ends our proof of Theorem 1.2(ii).

(iii) Suppose that any positive integer $n\eq9\pmod{20}$ can be written as $5x^2+5y^2+z^2$
with $x,y,z\in\Z$ and $2\nmid z$. Below we prove by induction that any $n\in\Z^+$ can be written as
$x^2+y^2+z^2+w^2$ with $x,y,z,w\in\Z$ and $x+3y\in\{4^k:\ k\in\N\}$.

It is easy to verify that each $n=1,\ldots,15$ can be written as $x^2+y^2+z^2+w^2$ with $x,y,z,w\in\Z$ and $x+3y\in\{4^k:\ k\in\N\}$.

Now let $n\in\Z^+$ with $n\gs16$, and assume that each $m=1,\ldots,n-1$ can be written as $x^2+y^2+z^2+w^2$ with $x,y,z,w\in\Z$ and $x+3y\in\{4^k:\ k\in\N\}$. If $16\mid n$, then by the induction hypothesis there are $x,y,z,w\in\Z$
such that $n/16=x^2+y^2+z^2+w^2$ and $x+3y\in\{4^k:\ k\in\N\}$, and hence
$n=(4x)^2+(4y)^2+(4z)^2+(4w)^2$ with $4x+3(4y)=4(x+3y)\in\{4^k:\ k\in\N\}$.

Below we suppose that $16\nmid n$.

{\it Case} 1. $2\nmid n$.

In this case, $10n-1\eq 9\pmod{20}$ and hence $10n-1=5u^2+5v^2+x^2$ for some $u,v,x\in\Z$ with $2\nmid x$.
As $u^2+v^2$ is even, both $y=(u+v)/2$ and $z=(u-v)/2$ are integers. Note that $10n-1=x^2+10y^2+10z^2$.

{\it Case} 2. $n=2m$ for some $m\in\Z^+$.

Note that $8\nmid m$ since $16\nmid n$. For $m=9,10,11,12$ we can easily verify the desired result.
Assume that $m\gs13$. Then $5m-64\gs 5\times 13-64>0$.
If $m$ is odd, then either $5m-4$ or $5m-64$ is not congruent to 7 mod 8.
If $2\|m$, then $5m-4\eq 5m-64\eq2\pmod4$.
If $4\|m$ and $(5m-64)/4=5m/4-16\eq7\pmod 8$, then $(5m-4)/4=5m/4-1\eq 6\pmod 8$.
So, for some $\da\in\{2,4\}$ we have $5m-(\da^2/2)^2\not\in E_0$, and hence by (2.3) there are $u,v,w\in\Z$
such that $5m-\da^4/4=5u^2+5v^2+w^2$.
Clearly,
$$10n-\da^4=4\l(5m-\f{\da^4}4\r)=10(2u^2+2v^2)+(2w)^2=(2w)^2+10(u+v)^2+10(u-v)^2.$$

In view of the above, for some $\da\in\{1,2,4\}$ we have $10n-\da^4=x^2+10y^2+10z^2$ for some $x,y,z\in\Z$.
As $x\eq \pm3\da^2\pmod{10}$, we may simply assume that $x=10w+3\da^2$ for some $w\in\Z$. Thus
$$10n-\da^4=(10w+3\da^2)^2+10y^2+10z^2$$
and hence
$$n=10w^2+6\da^2w+\da^4+y^2+z^2=(3w+\da^2)^2+(-w)^2+y^2+z^2$$
with $(3w+\da^2)+3(-w)=\da^2\in\{4^k:\ k\in\N\}$. This concludes the induction proof of Theorem 1.1(iii).
 \qed

\proclaim{Lemma 2.5} We have
$$E(1,3,6)=\{3q+2:\ q\in\N\}\cup\{4^k(16l+14):\ k,l\in\N\}\tag2.5$$
and
$$E(2,3,6)=\{3q+1:\ q\in\N\}\cup E_0.\tag2.6$$
\endproclaim
\Remark\ 2.3. (2.5) and (2.6) are known results, see, e.g., L. Dickson [D39, pp.\,112-113].

\medskip
\noindent{\it Proof of Theorem 1.3}. (i) We can easily verify the desired result for $n=1,\ldots,15$.

Now let $n\gs 16$ and assume that each $m=1,2,\ldots,n-1$ can be written as $x^2+y^2+z^2+w^2$ $(x,y,z,w\in\N$)
with $|x+y-z|\in\{4^k:\ k\in\N\}$.

We first suppose that $16\mid n$. By the induction hypothesis, there are $x,y,z,w,k\in\N$ for which $n/16=x^2+y^2+z^2+w^2$
with $|x+y-z|=4^k$, and hence $n=(4x)^2+(4y)^2+(4z)^2+(4w)^2$ with $|4x+4y-4z|=4^{k+1}$.

Now we suppose that $16\nmid n$. If $16<n<86$ then we can verify the desired result via a computer.
Thus we simply let $n\gs 86$ and hence $3n\gs 258>16^2$.
Let $\da=0$ if $3n-1\not\in E_0$.
In the case $3n-1\in E_0$, we let $\da=1$ if $n-6$ is not an odd square, and $\da=2$ otherwise.
By Lemmas 2.2 and 2.5, if $3n-1\in E_0$ then $3n-16,3n-16^2\not\in E(2,3,6)$.
As $3n>16^\da$ and $3n-16^{\da}\not\in E(2,3,6)$,
there are $x,y\in\N$ and $z\in\Z$ such that
$$3n-16^\da=3x^2+6y^2+2(3z-4^\da)^2=3(x^2+2y^2+2(3z^2-2\times4^{\da}z))+2\times16^{\da}$$
and hence
$$n=x^2+2y^2+6z^2-4^{\da+1}z+16^\da=x^2+(y+z)^2+(z-y)^2+(4^\da-2z)^2$$
with $(y+z)+(z-y)+(4^\da-2z)=4^{\da}$.

When $z\gs 2^{2\da-1}$, we have $2z\gs 4^{\da}$, hence
$$(y+z)+|z-y|-(2z-4^{\da})=4^{\da}\ \ \t{if}\ y\ls z,$$
and
$$||z-y|+(2z-4^{\da})-(y+z)|=4^{\da}\ \ \t{if}\ y>z.$$
When $y\gs z$ and $z<2^{2\da-1}$,
$$|y+z|+(4^\da-2z)-(y-z) = 4^\da\ \ \t{if}\ y+z\gs0,$$
and
$$|(y-z)+|y+z|-(4^\da-2z)|=4^\da\ \ \t{if}\ y+z<0.$$

Below we assume $0\ls y<z<2^{2\da-1}$.
Clearly $\da>0$. If $\da=1$, then we must have $y=0$ and $z=1$, hence
$$n=x^2+(y+z)^2+(z-y)^2+(4^\da-2z)^2=x^2+1^2+1^2+(4-2)^2=x^2+6.$$
If $2\mid x$, then $3(x^2+6)-1\eq 3\times 6-1\eq1\pmod4$ and hence $3(x^2+6)-1\not\in E_0$.
Thus $\da\not=1$ by the definition of $\da$. So $\da=2$ and $0\ls y<z<8$. As $n-6=x_0^2$ for a positive odd integer, we have
$$n=x^2+(y+z)^2+(z-y)^2+(16-2z)^2\eq7\pmod8,$$
hence $2\nmid x$, $2\nmid y\pm z$ and $2\nmid z$.

If $z=1$, then $y=0$ and
$$n=x^2+1^2+1^2+14^2=x^2+13^2+5^2+2^2$$
with $13+5-2=4^2$. If $z$ is $3$ or $5$, then
$$n=x^2+(y+z)^2+(z-y)^2+(16-2z)^2$$
with $|y+z+(z-y)-(16-2z)|=|4(z-4)|=4$.

Now we handle the remaining case $z=7$. Note that $y\in\{0,2,4,6\}$.
If $y=2$, then
$$n=x^2+(2+7)^2+(7-2)^2+(16-14)^2=x^2+6^2+5^2+7^2$$
with $6+5-7=4$.
If $y=4$, then
$$n=x^2+(4+7)^2+(7-4)^2+(16-14)^2=x^2+9^2+2^2+7^2$$
with $9+2-7=4$. If $y=6$, then
$$n=x^2+(6+7)^2+(7-6)^2+(16-14)^2=x^2+11^2+7^2+2^2$$
with $11+7-2=4^2$.

In the case $y=0$ and $z=7$, we have
$$x_0^2+6=n=x^2+(0+7)^2+(7-0)^2+(16-14)^2=x^2+102$$
and hence
$$\f{x_0-x}2\cdot\f{x_0+x}2=\f{102-6}4=24.$$
As $x_0$ and $x$ are positive and odd, $(x_0-x)/2\not\eq (x_0+x)/2\pmod 2$, hence either
$$\f{x_0-x}2=1,\ \f{x_0+x}2=24, \ \t{and thus}\ n=x_0^2+6=25^2+6,$$
or
$$\f{x_0-x}2=3,\ \f{x_0+x}2=8, \ \t{and thus}\ n=x_0^2+6=11^2+6.$$
Note that
$$11^2+6=1^2+5^2+10^2+1^2\ \t{with}\ |1+5-10|=4$$
and
$$25^2+6=1^2+5^2+22^2+11^2\ \t{with}\ |1+5-22|=4^2.$$

In view of the above, we have completed the induction proof of Theorem 1.3(i).
\medskip

(ii) Via a computer we can easily verify the desired result for every $n=0,1,\ldots,63$.

Now fix an integer $n\gs 64$ and assume that the desired result holds for all smaller values of $n$.

If $64\mid n$, then by the induction hypothesis we can write $n/64$ as $x^2+y^2+z^2+w^2$
with $x,y,z,w\in\N$ and $x+y-z\in\{\pm a8^k:\ k\in\N\}\cup\{0\}$, and hence
$$n=(8x)^2+(8y)^2+(8z)^2+(8w)^2\ \ \t{with}\ 8x+8y-8z\in\{\pm a8^k:\ k\in\N\}\cup\{0\}.$$

Below we suppose that $64\nmid n$.

By Lemma 2.4, $3n-(a\da)^2\not\in E_0$ for some $\da\in\{0,1,8\}$ satisfying (2.4). In view of (2.6),
for some $x,y,z\in\Z$ we have $3n-(a\da)^2=3x^2+6y^2+2(3z-a\da)^2$ and hence
$$n=x^2+(y+z)^2+(z-y)^2+(2z-a\da)^2.$$
When $z\gs a\da/2$, obviously $y+z\gs0$, $2z-a\da\gs0$,
$$(y+z)+(z-y)-(2z-a\da)=a\da\in\{\pm a8^k:\ k\in\N\}\cup\{0\}\ \ \t{if}\ y\ls z,$$
and
$$(y-z)+(2z-a\da)-(y+z)=-a\da\in\{\pm a8^k:\ k\in\N\}\cup\{0\}\ \ \t{if}\ y>z.$$
If $z<a\da/2$ and $y\gs |z|$, then $\{a\da-2z,y+z,y-z\}\se\N$ and
$$(a\da-2z)+(y+z)-(y-z)=a\da\in\{\pm a 8^k:\ k\in\N\}\cup\{0\}.$$
If $z<a\da/2$ and $z\ls -y$, then $\{a\da-2z,y-z,-y-z\}\se\N$ and
$$(y-z)+(-y-z)-(a\da-2z)=-a\da\in\{\pm a8^k:\ k\in\N\}\cup\{0\}.$$

Now we consider the remaining case $y<z<a\da/2$. Since $a\da>2$, we have $\da=8$, and hence $a=1$ and $16\mid n$ by (2.4).
In the case $y<z<8/2=4$, the ordered triple $(y+z,z-y,8-2z)$
is among
$$(1,1,6),\ (2,2,4),\ (3,3,2),\ (3,1,4),\ (4,2,2),\ (5,1,2).$$
Note that $2+2-4=0=3+1-4$. If
$$(y+z,z-y,8-2z)\in\{(1,1,6),\, (3,3,2),\, (5,1,2)\},$$
then
$$n=x^2+(y+z)^2+(z-y)^2+(8-2z)^2\eq x^2+2\not\eq0\pmod4$$
which contradicts $16\mid n$.
This concludes the proof of Theorem 1.3(ii). \qed

\proclaim{Lemma 2.6} Suppose that $n\in\Z^+$ is the sum of three squares. Then $n=a^2+b^2+c^2$
for some $a,b,c\in\Z$ with $a+b\eq1\pmod3$.
\endproclaim
\Proof. Write $n=x^2+y^2+z^2$ with $x,y,z\in\Z$. If $x\eq y\eq z\eq0\pmod 3$, then $9\mid n$
and hence by [S16, Lemma 2.2(ii)] we can write $n$ as $\bar x^2+\bar y^2+\bar z^2$ with $\bar x,\bar y,\bar z\in\Z$ and $3\nmid \bar x\bar y\bar z$.
So, without loss of generality we may assume that $3\nmid x$.
If $x+y\eq0\pmod 3$, then $-x+y\not\eq 0\pmod 3$. Thus, we may simply suppose that $x+y\not\eq0\pmod 3$ and hence
$x+y$ is congruent to $1$ or $-1$ modulo $3$.
If $x+y\eq-1\pmod3$, then $n=(-x)^2+(-y)^2+z^2$ with $(-x)+(-y)\eq1\pmod 3$.
Therefore, there are $a,b,c\in\Z$ such that $n=a^2+b^2+c^2$ and $a+b\eq1\pmod3$. \qed

\medskip
\noindent{\it Proof of Theorem 1.4}.
(i) We prove the desired result by induction. For $n=1,2,\ldots,42$ we can verify the result directly via a computer.

Now let $n\gs 43$ and assume that any $m=1,\ldots,n-1$ can be written as $x^2+y^2+z^2+w^2$ $(x,y,z,w\in\Z$) with
$x+y+2z=c4^k$ for some $k\in\N$.

If $16\mid n$, then by the induction hypothesis there are $x,y,z,w\in\Z$ and $k\in\N$ such that $n/16=x^2+y^2+z^2+w^2$ and
$x+y+2z=c4^k$,  hence $n=(4x)^2+(4y)^2+(4z)^2+(4w)^2$ with $4x+4y+2(4z)=c4^{k+1}$.

Now suppose that $16\nmid n$.
\medskip

Clearly, $3n-2\not\eq 14\pmod{16}$. If $4\nmid 3n-2$, then $3n-2\not\in E(1,3,6)$ by (2.4).
If $4\mid 3n-2$, then $4\nmid 3n-32$. If $3n-32\eq14\pmod{16}$, then $3n-2\eq12\pmod{16}$ and hence $3n-2\not\in E(1,3,6)$.
Thus, for some $\da\in\{0,1\}$, we can write
$$3n-2\times16^\da=3x^2+6y^2+(3z-4^\da)^2=3x^2+6y^2+9z^2-6\times 4^\da z+16^\da$$
with $x,y,z\in\Z$. It follows that
$$n=x^2+2y^2+3z^2-2\times4^\da z+16^\da=x^2+(y+z)^2+(z-y)^2+(4^\da-z)^2$$
with $(y+z)+(z-y)+2(4^\da-z)=2\times 4^\da$. This proves the desired result with $c=2$.

Below we show the desired result for $c=1$.

{\it Case} 1. $4\nmid n$.

In this case, $6n-1\not\eq7\pmod 8$ and hence $6n-1\not\in E(2,3,6)$ by (2.6). So, for some $x,y,z\in\Z$ we have
$$6n-1=6x^2+3(2y+1)^2+2(3z+1)^2$$
and hence
$$n=x^2+2y^2+2y+3z^2+1=x^2+(y+z+1)^2+(z-y)^2+(-z)^2$$
with $(y+z+1)+(z-y)+2(-z)=1=4^0$.
\medskip

{\it Case} 2. $4\|n$.

In this case, $3n-8\eq 4\pmod 8$ and hence $3n-8\not\in E(1,3,6)$ by (2.5). Thus, for some $x,y,z\in\Z$ we have
$3n-8=3x^2+6y^2+(3z-2)^2$ and hence
$$n=x^2+2y^2+3z^2-4z+4=x^2+(y+z)^2+(z-y)^2+(2-z)^2$$
with $(y+z)+(z-y)+2(2-z)=4$.
\medskip

{\it Case} 3. $8\|n$.

Write $n=8q$ with $q$ odd. Note that $q\gs 6$ and $3n>128$ since $n\gs 43$.
By Lemma 2.2, for some $\da\in\{0,1\}$ we have $3q-16^\da\not\in E_0$
and hence
$$3n-8\times 16^\da=8(3q-16^\da)\not\in E(1,3,6)$$
by (2.5). Thus, for some $x,y,z\in\Z$ we have
$$3n-8\times 16^\da=3x^2+6y^2+(3z-2\times 4^{\da})^2$$
and hence
$$n=x^2+2y^2+3z^2-4^{\da+1}z+4\times16^{\da}=x^2+(y+z)^2+(z-y)^2+(2\times4^\da-z)^2$$
with $(y+z)+(z-y)+2(2\times 4^\da-z)=4^{\da+1}$.

The induction proof of Theorem 1.4(i) is now completed.
\medskip

(ii) Let $m\in\{2,3\}$. For $n=1,2,\ldots,4^m-1$ we can easily verify that $n^{m-1}$ can be written as $x^2+y^2+z^2+w^2$
with $x,y,z,w\in\Z$ and $x+2y+2z\in\{2^{km}:\ k\in\N\}$.

Now let $n\gs 4^m$ and assume that for each $n_0=1,2,\ldots,n-1$ we can write $n_0^{m-1}$ as $x^2+y^2+z^2+w^2$ $(x,y,z,w\in\Z$)
with $x+2y+2z$ a power of $2^m$.

If $4^m\mid n^{m-1}$, then by the induction hypothesis there are $x,y,z,w\in\Z$ and $k\in\N$ for which $n^{m-1}/4^m=x^2+y^2+z^2+w^2$
with $x+2y+2z=2^{km}$, and hence $n^{m-1}=(2^mx)^2+(2^my)^2+(2^mz)^2+(2^mw)^2$ with $2^mx+2(2^my)+2(2^mz)=2^{(k+1)m}$.

Now we suppose that $4^m\nmid n^{m-1}$. By Lemmas 2.1-2.2, if $16\nmid n^{m-1}$, then for some $\da\in\{0,1\}$  we can write $9n^{m-1}-4^{\da m}$ as $a^2+b^2+c^2$ with $a,b,c\in\Z$. When $16\mid n^{m-1}$, we must have $m=3$ and $4\|n$,
thus by taking $\da=1$ we find that $9n^{m-1}-4^{\da m}=4^2(9(n/4)^2-4)$ can be written as $a^2+b^2+c^2$ with $a,b,c\in\Z$.
Clearly, we cannot have $3\nmid abc$.
Without loss of generality, we assume that $c=3w$ with $w\in\Z$. As $a^2+b^2\eq-16^\da\eq2\pmod 3$, we must have $3\nmid ab$.
We may simply suppose that $a=3u+2^{\da m+1}$ and $b=3v-2^{\da m+1}$ with $u,v\in\Z$. (Note that if $x\eq1\eq-2\pmod3$ then $-x\eq2\pmod 3$.)
Since
$$\align &12\times 2^{\da m}u-12\times 2^{\da m} v+8(2^{\da m})^2
\\\eq&(3u+2\times2^{\da m})^2+(3v-2\times2^{\da m})^2=a^2+b^2\eq-(2^{\da m})^2\pmod 9,
\endalign$$
we must have $u\eq v\pmod 3$.
Set
$$y=\f{2u+v}3\ \ \t{and}\ \ z=\f{u+2v}3.$$
Then
$$\align 9n^{m-1}-4^{\da m}=&a^2+b^2+c^2=(3u+2^{\da m+1})^2+(3v-2^{\da m+1})^2+9w^2
\\=&(3(2y-z)+2^{\da m+1})^2+(3(2z-y)-2^{\da m+1})^2+9w^2
\\=&9(2y-z)^2+9(2z-y)^2+3\times 2^{\da m+2}((2y-z)-(2z-y))
\\&+8\times 4^{\da m}+9w^2
\endalign$$
and hence
$$\align n^{m-1}=&(2y-z)^2+(2z-y)^2+2^{\da m+2}(y-z)+w^2+4^{\da m}
\\=&(2y-2z+2^{\da m})^2+(-y)^2+z^2+w^2
\endalign$$
with $(2y-2z+2^{\da m})+2(-y)+2z=2^{\da m}$. This proves Theorem 1.4(ii).

(iii) For any $k\in\N$, we obviously have
$$2^{2k+1}=(2^{k})^2+(2^{k})^2+0^2+0^2\ \t{with}\ 2^{k}+2\times 2^{k}+2\times0=3\times 2^{k},$$
and
$$2^{2k+2}=(-2^{k})^2+(2^{k})^2+(2^{k})^2+(2^{k})^2\ \t{with}\ (-2^{k})+2\times 2^{k}+2\times2^{k}=3\times2^{k}.$$
So, it suffices to prove the first assertion in Theorem 1.4(iii).

For $n\in\{1,\ldots,15\}$ with $n\not=1,8$, we can easily verify that $n$ can be written as $x^2+y^2+z^2+w^2$ $(x,y,z,w\in\Z$)
with $x+2y+2z\in\{3\times 4^k:\ k\in\N\}$. For example, $4=(-1)^2+1^2+1^2+1^2$ with $(-1)+2\times1+2\times1=3$.

Now let $n\gs 16$ with $n\not\in S=\bigcup_{k\in\N}\{2^{4k},2^{4k+3}\}$,
and assume that each $m=1,2,\ldots,n-1$ with $m\not\in S$ can be written as $x^2+y^2+z^2+w^2$ $(x,y,z,w\in\Z$)
with $x+2y+2z\in\{3\times 4^k:\ k\in\N\}$.

If $16\mid n$, then by the induction hypothesis there are $x,y,z,w\in\Z$ and $k\in\N$ for which $n/16=x^2+y^2+z^2+w^2$
with $x+2y+2z=3\times 4^k$, and hence $n=(4x)^2+(4y)^2+(4z)^2+(4w)^2$ with $4x+2(4y)+2(4z)=3\times4^{k+1}$.

Now we suppose that $16\nmid n$. By Lemmas 2.1-2.2, for some $\da\in\{0,1\}$  we can write $n-16^\da$ as the sum of three squares.
Combining this with Lemma 2.6, we see that $n-16^\da=a^2+b^2+c^2$ for some $a,b,c\in\Z$ with $a+b\eq1\pmod 3$.
Let $u=a-2^{2\da+1}$ and $v=b-2^{2\da+1}$.
Then $u+v\eq a-2+(b-2)\eq 0\pmod 3$, and
$$y=\f{2u-v}3\ \ \t{and}\ \ z=\f{2v-u}3$$
are both integers. Observe that
$$n-16^\da=(u+2\times 4^\da)^2+(v+2\times 4^\da)^2+c^2=(2y+z+2\times 4^\da)^2+(y+2z+2\times 4^\da)^2+c^2$$
and hence
$$n=(2y+2z+3\times4^\da)^2+(-y)^2+(-z)^2+c^2$$
with $(2y+2z+3\times4^\da)+2(-y)+2(-z)=3\times 4^\da$.

\medskip

So far we have completed the proof of Theorem 1.4. \qed

\heading{3. Proofs of Theorems 1.5 and 1.6}\endheading

Let us first recall some known results on ternary quadratic forms.

\proclaim{Lemma 3.1} We have
$$\align E(1,2,4)=&\{4^k(16l+14):\ k,l\in\N\},\tag3.1
\\E(1,6,9)=&\{3q+2:\ q\in\N\}\cup\{9^k(9l+3):\ k,l\in\N\}.\tag3.2
\\E(2,3,12)=&\{16q+6:\ q\in\N\}\cup\{9^k(3l+1):\ k,l\in\N\},\tag3.3
\\E(1,5,10)=&\bigcup_{k,l\in\N}\{25^k(5l+2),25^k(5l+3)\},\tag3.4
\\E(2,5,10)=&\{8q+3:\ q\in\N\}\cup\bigcup_{k,l\in\N}\{25^k(5l+1),25^k(5l+4)\}.\tag3.5
\endalign$$
\endproclaim
\Remark\ 3.1. (3.1)-(3.5) can be found in L. E. Dickson [D39, pp.\,112-113].
\medskip

\proclaim{Lemma 3.2} Let $n\in\Z^+$ and $\da\in\{0,1\}$. Then $6n+1=x^2+3y^2+6z^2$ for some $x,y,z\in\Z$ with $x\eq\da\pmod2$.
\endproclaim
\Remark\ 3.2. This appeared in Sun [S17a, Remark 3.1].
\medskip

We also need the following result in [S15].

\proclaim{Lemma 3.3} {\rm (i)} Any $n\in\N$ with $n\eq4\pmod{12}$ can be written as $x^2+3y^2+3z^2$ with $x,y,z\in\Z$ and $2\nmid x$.

{\rm (ii)} For $n\in\N$ with $n\eq4\pmod 8$, we have
$$|\{(x,y)\in\Z^2:\ x^2+3y^2=n\ \t{and}\ 2\nmid xy\}|=\f23|\{(x,y)\in\Z^2:\ x^2+3y^2=n\}|.\tag3.6$$
\endproclaim
\Remark\ 3.3. For parts (i) and (ii) one may consult Theorem 1.7(iii) and Lemma 3.2 of Sun [S15].

\proclaim{Lemma 3.4} Let $n\in\Z^+$ with $n\eq1,2,4\pmod7$ and $n\not\eq2\pmod3$. Then $n=x^2+7y^2+14z^2$ for some $x,y,z\in\Z$.
\endproclaim
\Proof. In light of [BIS], the genus $G$ (of discriminant $-392$) containing the class of the form $x^2+7y^2+14z^2$ contains only one other class, namely the class containing the form $2x^2+7y^2+7z^2$. By Jones [J31a, p.\,99]
and [J31b, p.\,123], a positive integer is represented by a form in $G$ if and only if it does not belong to the set
$$\l\{7^{2k}(7m+r):\ k,m\in\N\ \t{and}\ r\in\{3,5,6\}\r\}.$$
As $n\not\eq3,5,6\pmod 7$, it is represented by a form in $G$, i.e., $n$ is represented by $x^2+7y^2+14z^2$ or $2x^2+7y^2+7z^2$ or both.

Suppose that $n=2x^2+7y^2+7z^2$ for some $x,y,z\in\Z$. If $y^2,z^2\not\eq x^2\pmod 3$, then $y^2\eq z^2\pmod3$
and hence $2x^2+7y^2+7z^2\eq2x^2+2y^2\eq2\pmod3$. As $n\not\eq2\pmod3$, $x^2$ is congruent to $y^2$ or $z^2$ modulo $3$.
Without loss of generality we may simply assume that $x\eq z\pmod3$. (If $x\eq-z\pmod3$ then we may replace $z$ by $-z$.)
Thus both $u=(2x+7z)/3$ and $v=(x-z)/3$ are integers. Note that $u^2+7y^2+14v^2=2x^2+7y^2+7z^2=n$.

In view of the above, we immediately obtain the desired result. \qed
\medskip

\Remark\ 3.4. I. Kaplansky [K] reported that $2$, $74$ and $506$
are the only positive integers $n\ls 10^5$ with $n\eq1,2,4\pmod 7$ which cannot be represented by $x^2+7y^2+14z^2$ with $x,y,z\in\Z$. We guess that any positive integer $n\eq1\pmod7$ can be written as $x^2+7y^2+14z^2$ with $x,y,z\in\Z$.

\medskip
\noindent {\it Proof of Theorem} 1.5. (i) By (3.1), for each $n\in\Z^+$ we can write $2n-1=x^2+2z^2+4w^2$ with $x,z,w\in\N$. As $x$ is odd, we may write $x=2y+1$ with $y\in\N$.
Thus $2n-1=(2y+1)^2+2z^2+4w^2$ and hence $n=x^2+y^2+z^2+2w^2$ with $x=y+1$.

By the above, we can write any $n\in\Z^+$ as $x^2+y^2+z^2+2w^2$ $(x,y,z,w\in\Z$) with $x+y=1$. Clearly, $z$ is congruent to $x$ or $y$ modulo $2$.
Without loss of generality, we may assume that $y\eq z\pmod 2$. Then
$$n=x^2+2\l(\f{y+z}2\r)^2+2\l(\f{y-z}2\r)^2+2w^2\ \ \t{with}\ x+\f{y+z}2+\f{y-z}2=1.$$

(ii) Let $n\in\Z^+$.  By (3.3), $6n-1\not\in E(2,3,12)$. So $6n-1=2u^2+3v^2+12w^2$ for some $u,v,w\in\Z$. As $2\nmid v$, we may write $v=2y+1$ with $y\in\Z$.
Since $3\nmid u$, we may write $u$ or $-u$ as $3z+1$ with $z\in\Z$. Thus
$$6n-1=2(3z+1)^2+3(2y+1)^2+12w^2$$
and hence
$$n=(y+z+1)^2+(z-y)^2+(-z)^2+2w^2$$
with $(y+z+1)+(z-y)+2(-z)=1$.

By the above, $n=x^2+y^2+z^2+2w^2$ for some $x,y,z,w\in\Z$ with $x+y+2z=1$. As $x+y$ is odd, $x$ or $y$ is even.
Thus $n=u^2+(2v)^2+z^2+2w^2$ for some $u,v\in\Z$ with $u+2v+2z=1$.
Clearly, $x$ or $y$ is congruent to $z$ modulo $2$.
Without any loss of generality, we may assume that $y\eq z\pmod 2$. Observe that
$$n=x^2+2\l(\f{z-y}2\r)^2+2\l(\f{y+z}2\r)^2+2w^2$$
with
$$ x+\f{z-y}2+3\f{y+z}2=x+y+2z=1.$$

(iii) As $4n-1\not\in E(1,2,4)$ by (3.1), we have $4n-1=u^2+2v^2+4x^2$ for some $u,v,x\in\Z$.
Since $u$ or $-u$ is congruent to $1$ modulo 4, without loss of generality we may assume that $u=4y+1$ with $y\in\Z$.
As $u^2\not\eq-1\pmod 4$, we have $v=2z+1$ for some $z\in\Z$. Thus
$4n-1=4x^2+(4y+1)^2+2(2z+1)^2$ and hence
$$n=x^2+4y^2+2y+2z^2+2z+1=x^2+(y+z+1)^2+(y-z)^2+2(-y)^2$$
with $(y+z+1)+(y-z)+2(-y)=1$.

By the above, there are $x,y,z,w\in\Z$ with $n=x^2+y^2+z^2+2w^2$ and $y+z+2w=1$.
As $y+z$ is odd, one of $y$ and $z$ is odd and another is even. If $z=2t$ with $t\in\Z$, then $n=x^2+y^2+4t^2+2w^2$ with $y+2t+2w=1$.
Clearly $x$ is congruent to $y$ or $z$ modulo 2; if $x\eq y\pmod 2$ then
$$n=2\l(\f{x+y}2\r)^2+2\l(\f{y-x}2\r)^2+z^2+2w^2$$
with
$$\f{x+y}2+\f{y-x}2+z+2w=y+z+2w=1.$$
This proves Theorem 1.5(iii).

(iv) We first suppose that $n\not\eq1\pmod 8$. Then $5n-2\not\eq 3\pmod 8$ and hence $5n-2\not\in E(2,5,10)$ by (3.5). So
$5n-2=2v^2+5x^2+10z^2$ for some $v,x,z\in\Z$. As $v^2\eq -1\eq 2^2\pmod 5$, $v$ or $-v$ is congruent to $-2$ modulo $5$.
Without any loss of generality, we may assume that $v=5y-2$ with $y\in\Z$. Thus
$$5n=2+2(5y-2)^2+5x^2+10z^2$$
and hence
$$n=x^2+10y^2-8y+2+2z^2=x^2+(y+z)^2+(y-z)^2+2(1-2y)^2$$
with $(y+z)+(y-z)+(1-2y)=1$.

Now we handle the case $n\eq1\pmod 8$. As $10n-4\not\in E(1,5,10)$ by (3.4), there are $u,v,x\in\Z$
with $10n-4=u^2+5v^2+10x^2$. If $u$ and $v$ are both even, then $10x^2\eq 10n-4\eq 2\pmod 4$ and hence $2\nmid x$,
thus
$$\l(\f u2\r)^2+\l(\f v2\r)^2\eq\l(\f u2\r)^2+5\l(\f v2\r)^2=\f{10(n-x^2)-4}4\eq-1\pmod 4$$
which is impossible. Thus $2\nmid uv$ and we can write $v=2z+1$ with $z\in\Z$. Since $u^2\eq1\pmod 5$,
without loss of generality we may assume that $u=10y+1$ for some $y\in\Z$. Thus
$10n=4+(10y+1)^2+5(2z+1)^2+10x^2$ and hence
$$n=x^2+10y^2+2y+2z^2+2z+1=x^2+(y+z+1)^2+(y-z)^2+2(-2y)^2$$
with $(y+z+1)+(y-z)+(-2y)=1$.

(v) Let $n\in\Z^+$. By Lemma 3.2, we can write $6n-5=u^2+3v^2+6x^2$ with $u,v,x\in\Z$ and $2\nmid u$.
As $u$ or $-u$ is congruent to $-1$ modulo $6$, we may simply assume that $u=6w-1$ for some $w\in\Z$.
Since $v$ is even, $v=2z$ for some $z\in\Z$. Thus
$$6n-5=(6w-1)^2+3(2z)^2+6x^2=36w^2-12w+1+12z^2+6x^2$$
and hence
$$n=6w^2-2w+1+2z^2+x^2=x^2+(1-w)^2+2z^2+5w^2$$
with $(1-w)+w=1$.

Suppose that $n\in\Z^+$ and $n\not\eq2\pmod3$. Let $\da$ be $1$ or $2$. By Lemma 3.4, we have
$7n-6/\da=v^2+7x^2+14z^2$ for some $v,x,z\in\Z$. Since $v^2\eq\da^2\pmod 7$, without loss of generality we may assume that
$v=7w-\da$ for some $w\in\Z$. Thus
$$\align n=&\f{7x^2+14z^2+(7w-\da)^2+6/\da}7=x^2+2z^2+7w^2-2\da w+\f17\l(\da^2+\f 6{\da}\r)
\\=&x^2+(1-\da w)^2+2z^2+\f{6}{\da}w^2
\endalign$$
with $(1-\da w)+\da w=1$.

(vi) As $3n-2\not\in E(1,6,9)$ by (3.2), there are $u,v,w\in\Z$ such that
$$3n-2=(3u+2)^2+6(v+1)^2+9w^2$$
and hence
$$n=3u^2+4u+4+2v^2+4v+3w^2=(u+v+2)^2+(u-v)^2+(-u)^2+3w^2$$
with $(u+v+2)+(u-v)+2(-u)=2$.

(vii) Note that $5=1^2+(-1)^2+1^2+2\times1^2$ with $1+(-1)+1=1^2$.

Below we let $n\in\N$ with $n\gs 6$. Choose $c\in\{1,4\}$ such that $6n-2c^2\eq0\pmod 4$.
Note that $6n-2c^2\eq c^2\eq1\pmod 3$ and hence $6n-2c^2\eq4\pmod{12}$. By Lemma 3.3(i), there are $r,s,t\in\Z$ with $r$ odd such that $6n-2c^2=r^2+3s^2+3t^2$.
As $s\not\eq t\pmod 2$, without loss of generality we simply assume that $2\nmid s$ and $t=2w$ with $w\in\Z$. Since $r^2+3s^2\eq1+3=4\pmod 8$, by Lemma 3.3(ii) we can
write $r^2+3s^2=u^2+3v^2$ with $u,v\in\Z$ and $u\eq v\eq c\pmod 2$. Clearly $v=2y+c$ for some $y\in\Z$.
As $u$ or $-u$ is congruent to $c$ modulo 3, we may write $u$ or $-u$ as $6z+c$ with $z\in\Z$. Thus
$$6n-2c^2=(6z+c)^2+3(2y+c)^2+3(2w)^2=12y^2+12c(y+z)+36z^2+12w^2+4c^2$$
and hence
$$n=2y^2+2z^2+2c(y+z)+c^2+4z^2+2w^2=(y+z+c)^2+(z-y)^2+(-2z)^2+2w^2$$
with $(y+z+c)+(z-y)+(-2z)=c\in\{t^2:\ t=1,2\}$.

(viii) For $n=8,9,10$ we can easily verify that $n$ can be written as $x^2+y^2+z^2+2w^2$ $(x,y,z,w\in\Z$) with $x+y+2z\in\{2t^2:\ t=1,2\}$.

Now assume that $n\gs11$. As $6n-2^2\not\eq 6n-8^2\pmod{16}$, by (3.3) we have $6n-c^2\not\in E(2,3,12)$ for some $c\in\{2t^2:\ t=1,2\}$.
Then we can write $6n-c^2=2(3z+c)^2+3(2y+c)^2+12w^2$ with $y,z,w\in\Z$. It follows that
$$n=(y+z+c)^2+(z-y)^2+(-z)^2+2w^2$$
with $(y+z+c)+(z-y)+2(-z)=c\in\{2t^2:\ t=1,2\}$.

\medskip
Combining the above, we have completed the proof of Theorem 1.5. \qed

\proclaim{Lemma 3.5} We have the new identity
$$\aligned &(a^2+ab+b^2)(av^2+b(s^2+t^2+u^2))
\\=&a(b(s+t-u)+(a-b)v)^2+b(bs+au+av)^2
\\&+b(at+bu-av)^2+b(as-bt-av)^2.
\endaligned\tag3.7$$
In particular,
$$7(s^2+t^2+u^2+2v^2)=x^2+y^2+z^2+2w^2,\tag3.8$$
where
$$x=s+2u+2v,\ y=-2t-u+2v,\ z= 2s-t-2v,\ w=s+t-u+v.\tag3.9$$
\endproclaim
\Proof. By expanding and simplifying the right-hand side of (3.7), we see that (3.7) does hold.
Putting $a=2$ and $b=1$ in (3.7), we immediately obtain (3.8). \qed

\proclaim{Lemma 3.6} Let $m\in\{2,3\}$. For any integer $n\gs2\times7^{2m-1}$, there are integers $x,y,z$ and $w\in\{0,7^m\}$ such that
$$7n=x^2+y^2+z^2+2w^2\ \ \t{with}\ x+3y\eq0\pmod 7.\tag3.10$$
\endproclaim
\Proof. Note that $7n\gs 2\times 7^{2m}=2\times 49^m$. If $7n\in E(1,1,1)$, then $7n=4^k(8l+7)$ for some $k,l\in\N$,
and hence $7n-2\times49^m\not\in E(1,1,1)$. So, for some $w\in\{0,7^m\}$ we can write $7n-2w^2$ as $x^2+y^2+z^2$
with $x,y,z\in\Z$. Note that $x^2+y^2+z^2\eq0\pmod 7$.

Since $x^2+y^2\eq 6z^2\pmod 7$, we have
$x^2\eq (2z)^2$ and $y^2\eq(3z)^2$, or  $x^2\eq (3z)^2$ and $y^2\eq(2z)^2$.
Without any loss of generality, we may assume that $x\eq2z\pmod 7$ and $y\eq-3z\pmod 7$,
hence $x+3y\eq 2z-9z\eq0\pmod7$. This concludes our proof. \qed

\medskip
\noindent{\it Proof of Theorem 1.6}. For $n=0,1,\ldots,2\times 7^{2m-1}-1$ we can verify the desired results via a computer.

Below we fix an integer $n\gs 2\times 7^{2m-1}$.

(i) By Lemma 3.6, there are $x,y,z\in\Z$ and $w\in\{0,7^m\}$ satisfying (3.10).
Note that $x+3y\eq w\eq0\pmod 7$ and
$$\align z^2=7n-x^2-y^2-2w^2\eq& -x^2-y^2
\\\eq& 6x^2+27y^2=(2x-3y)^2+2(x+3y)^2
\\\eq&(2x-3y)^2\pmod 7.
\endalign$$
Without loss of generality, we may assume that $z\eq2x-3y\pmod 7$.
Define
$$\cases s=\f{x+2z+2w}{7},\\t=\f{2w-2y-z}7,\\u=\f{2x-y-2w}7,\\v=\f{x+y-z+w}7.\endcases\tag3.11$$
It is easy to see that (3.9) holds. Thus, by Lemma 3.5 we have (3.8) and hence
$$n=\f{7n}7=\f{x^2+y^2+z^2+2w^2}{7}=s^2+t^2+u^2+2v^2.$$
As $x\eq-3y\pmod7$ and $z\eq 2x-3y\pmod7$, we see that $s,t,u,v\in\Z$ and that $s+t-u+v=w$ is an $m$-th power.

(ii) Choose $x\in\{0,7^m\}$ such that $7n-x^2$ is odd. By Dickson [D39, pp.\,112-113],
$$E(1,1,2)=\{4^k(16l+14):\ k,l\in\N\}.$$
So $7n-x^2=y^2+z^2+2w^2$ for some $y,z,w\in\Z$. As $y^2+z^2\eq 5w^2\pmod 7$, we have $y^2\eq(2w)^2\pmod 7$ and $z^2\eq w^2\pmod 7$,
or $y^2\eq w^2\pmod7$ and $z^2\eq(2w)^2\pmod7$. Without any loss of generality, we may assume that $y\eq -2w\pmod 7$ and $z\eq -w\pmod 7$.
Now it is easy to see that the numbers $s,t,u,v$ given by (3.11) are all integral. Note that $n=s^2+t^2+u^2+2v^2$ by Lemma 3.5.
Clearly, $s+2u+2v=x$ is an $m$-th power.

The proof of Theorem 1.6 is now complete. \qed

\heading{4. Some open conjectures}\endheading

In this section, we pose 16 open conjectures on sums of squares for further research.

\proclaim{Conjecture 4.1} {\rm (i)} Any integer $n>1$ can be written as the sum of two squares, a power of three and a power of five; in other words, we have
$$\{a^2+b^2+3^c+5^d:\ a,b,c,d\in\N\}=\{2,3,4,\ldots\}.$$

{\rm (ii)} Each integer $n>1$ can be written as the sum of two squares and two central binomial coefficients; in other words, we have
$$\l\{a^2+b^2+\bi{2c}c+\bi{2d}d:\ a,b,c,d\in\N\r\}=\{2,3,4,\ldots\}.$$

{\rm (iii)} Any integer $n>5$ can be written as $a^2+b^2+2^c+5\times2^d$ with $a,b,c,d\in\N$.
\endproclaim
\Remark\ 4.1. See [S, A303656, A303540 and A303637] for related data, and note that $\bi{2k}k\sim 4^k/\sqrt{k\pi}$ as $k\to+\infty$.
We have verified parts (i)-(iii) for $n$ up to $2\times10^{10}$, $10^{10}$ and $5\times10^9$ respectively.
The author would like to offer 3500 US dollars as the prize for the first proof of part (i) of Conjecture 4.1.
In contrast with Conjecture 4.1(iii), R. Crocker [C] showed in 2008 that there are infinitely many positive integers not representable as the sum of two squares and at most two powers of 2 (see also [PT] for a simple proof). We also conjecture that any integer $n>1$ can be written as the sum of two triangular numbers and two powers of $5$ (cf. [S, A303389]).
\medskip

We also have the following conjecture on restricted sums of three squares.

\proclaim{Conjecture 4.2} {\rm (i)} Any $n\in\Z^+$ with $\ord_2(n)$ odd can be written as $x^2+y^2+z^2\ (x,y,z\in\Z)$
with $x+3y+5z$ a square (or twice a square).

{\rm (ii)} Any $n\in\N$ not of the form $4^k(8l+7)\ (k,l\in\N)$ can be written as $x^2+y^2+z^2$ with $x,y,z\in\Z$ such that $x+2y+3z$
is a square or twice a square.

{\rm (iii)} Let $n\in\N$. Then $8n+1$ can be written as $x^2+y^2+z^2$ with $x,y\in\Z$ and $z\in\Z^+$ such that $x+3y$ is a square.
Also, we can write $8n+6$ as $x^2+y^2+z^2$ with $x\in\Z$, $y,z\in\N$ and $2\nmid z$ such that $x+2y$ is a square.

{\rm (iv)} We can write any positive odd integer as $x^2+2y^2+3z^2$ with $x,y,z\in\Z$ such that $x+y+z$ is a square or twice a square.
\endproclaim
\Remark\ 4.2. See [S, A283269, A283273 and A283299] for related data.
It is known that any positive odd integer can be written as  $x^2+2y^2+3z^2$ with $x,y,z\in\Z$
(cf. [D39, pp.\,112-113]).
\medskip

As $2(x^2+y^2+z^2+w^2)=(x+y)^2+(x-y)^2+(z+w)^2+(z-w)^2$, Lagrange's four-square theorem
is equivalent to the fact that each positive odd integer can be written as the sum of four squares.
Our following conjecture provides some refinements of Lagrange's four-square theorem involving primes.

\proclaim{Conjecture 4.3} {\rm (i)} Any positive odd integer can be written as $x^2+y^2+z^2+w^2$
with $x,y,z,w\in\Z$ such that $p=x^2+3y^2+5z^2+7w^2$ and $p-2$ are twin prime.

{\rm (ii)} Any integer $n>1$ not divisible by $4$ can be written as $x^2+y^2+z^2+w^2$
($x,y,z,w\in\N$) such that $p=x+2y+5z$, $p-2$, $p+4$ and $p+10$ are all prime.

{\rm (iii)} Any positive odd integer can be written as $x^2+y^2+z^2+4w^2$
with $x,y,z,w\in\N$ such that $2^x+2^y+2^z+1$ is prime.

{\rm (iv)} Any odd integer $n>1$ can be written as $x^2+y^2+z^2+w^2$
($x,y,z,w\in\N$) such that $2^{x+y}+2^{z+w}+1$ is prime.
\endproclaim
\Remark\ 4.3. See [S, A290935, A291635, A291150 and A291191] for related data. For example,
$39 = 1^2 + 3^2 + 5^2 + 2^2$
with $1^2 + 3\cdot3^2 + 5\cdot5^2 + 7\cdot2^2 = 181$ and $181 - 2 = 179$ twin prime,
 $143=1^2 + 5^2 + 9^2 + 4\cdot3^2$ with $2^1 + 2^5 + 2^9 + 1 = 547$ prime, and
$$2\times6998538+1 = 122^2 + 220^2 + 208^2 + 3727^2$$
 with $2^{122+220} + 2^{208+3727} + 1 = 2^{342} + 2^{3935} + 1$ a prime of 1185 decimal digits.
 Clearly, part (i) of Conjecture 4.3 unifies Lagrange's four-square theorem and the twin prime conjecture.
\medskip

The following Conjectures 4.4-4.6 are mainly motivated by Theorems 1.2-1.4.

\proclaim{Conjecture 4.4} {\rm(i)} Any $n\in\Z^+$ can be written as $x^2+y^2+z^2+w^2\ (x,y,z,w\in\N)$ with $P(x,y,z)\in\{2^k:\ k\in\N\}$, whenever $P(x,y,z)$
is among the polynomials
$$\gather 2x-y,\ 2x-3y,\ x+(y-z)/3,\ 2x+(y-z)/3,\ 2x-2y-z,
\\ 4x-2y-z,\ 4x-3y-z,\ 4x-4y-3z,\ x+y-z,\ x+y-2z,\ x+2y-z,
\\ x+3y-z,\ x+3y-2z,\ x+3y-3z,\ x+3y-4z,\ x+3y-5z,
\\x+4y-z,\ x+4y-2z,\ x+4y-3z,\ x+4y-4z,\ x+5y-z,\ x+5y-2z,
\\ x+5y-4z,\ x+5y-5z,\ x+6y-3z,\ x+7y-4z,\ x+7y-7z,\ x+8y-z,
\\ x+9y-2z,\ 2x+3y-z,\ 2x+3y-3z,\ 2x+3y-4z,\ 2x+5y-z,
\\ 2x+5y-3z,\ 2x+5y-4z,\ 2x+5y-5z,\ 2x+7y-z,\ 2x+7y-3z,
\\ 2x+7y-7z,\ 2x+9y-3z,\ 2x+11y-5z,\ 3x+4y-3z,\ 7x+8y-7z.
\endgather$$

{\rm (ii)} Any $n\in\Z^+$ can be written as $x^2+y^2+z^2+w^2\ (x,y,z,w\in\N)$ with $Q(x,y,z,w)\in\{2^k:\ k\in\N\}$, whenever $Q(x,y,z,w)$
is among the polynomials
$$\gather x+y+2z-2w,\ x+y+2z-3w,\ x+y+2z-4w,\ x+2y+2z-3w,
\\ x+2y+3z-4w,\ x+2y+4z-3w,\ x+2y+6z-7w,\ x+3y+4z-4w,
\\ x+4y+6z-5w,\ 2x+3y+5z-4w,\ 2x+y-z-w,\ 2x+y-2z-w,
\\ 2x+y-3z-w,\ 2x+2y-3z-2w,\ 3x+y-3z-2w,\ 3x+y-4z-2w,
\\ 3x+2y-2z-w,\ 3x+2y-3z-w,\ 3x+2y-3z-2w,\ 3x+y-2z-w,
\\ 3x+2y-6z-w,\ 4x+y-2z-w,\ 4x+2y-2z-w,\ 4x+3y-2z-w,
\\4x+3y-3z-w,\ 4x+3y-4z-3w,\ 4x+3y-5z-w,\ 4x+3y-5z-2w,
\\ 5x+2y-2z-w,\ 5x+2y-3z-2w,\ 5x+2y-4z-3w,\ 5x+4y-5z-w,
\\ 7x+y-6z-2w,\ 8x+3y-3z-2w,\ 8x+3y-10z-w,\ 9x+y-4z-w.
\endgather$$
\endproclaim

\proclaim{Conjecture 4.5} {\rm (i)} Each $n\in\N$ can be written as $x^2+y^2+z^2+w^2$ with $x,y,z,w\in\N$
and $P(x,y,z,w)\in\{4^k:\ k\in\N\}\cup\{0\}$, whenever $P(x,y,z,w)$ is among the linear polynomials
$$\align &2x-y,\ x+y-z,\ x-y-z,\ x+y-2z,\ \ 2x+y-z,\ 2x-y-z,
\\&2x-2y-z,\ 2x+y-3z,\ 2x+2y-2z,\  2x+2y-4z,\ 3x-2y-z,
\\&x+3y-3z,\ 2x+3y-3z,\ 4x+2y-2z,\ 8x+2y-2z,
\\&2(x-y)+z-w,\ 4(x-y)+2(z-w).
\endalign$$

{\rm (ii)} Any $n\in\Z^+$ can be written as $x^2+y^2+z^2+w^2$ with $x,y,z,w\in\N$ and $|Q(x,y,z,w)|\in\{4^k:\ k\in\N\}$, whenever $Q(x,y,z,w)$
is among the polynomials
$$\align &x+y-3z,\ x+2y-3z,\ x+2y-4z,\ x+2y-5z,\ x+3y-3z,
\\&x+3y-5z,\ x+4y-2z,\ x+5y-2z,\ x+5y-7z,\ 2x+3y-4z,
\\& 2x+4y-6z,\ 2x+4y-10z,\ 2x+5y-4z,\ 4x+6y-6z,
\\&4x+6y-14z,\ x+4y-2z-w,\ x+8y-3z-w,
\\&2x+3y-3z-w,\ x+y+2z-2w,\ x+2y+3z-4w.
\endalign$$

{\rm (iii)} Each $n\in\N$ can be written as $x^2+y^2+z^2+w^2$ with $x,y,z,w\in\N$
and $R(x,y,z,w)\in\{\pm8^k:\ k\in\N\}\cup\{0\}$, whenever $R(x,y,z,w)$ is among the linear polynomials
$$\align &x+2y-2z,\ x+3y-3z,\ 2x+3y-3z,\ 4x+6y-6z,\ 4x+8y-8z,
\\&4x+2y-10z,\ 4x+12y-12z,\ 8x+12y-12z,\ 8x+4y-20z,\ 2x+y-z-w,
\\& 8x+4y-4z-4w,\ 2x+8y-4z-2w,\ 4x+y-2z-w,\ 4x+16y-8z-4w.
\endalign$$
\endproclaim

\proclaim{Conjecture 4.6} {\rm (i)} Let $a\in\Z^+$ and let $b,c,d\in\N$ with $a\gs b\gs c\gs d$ and $4\nmid \gcd(a,b,c,d)$. Then any $n\in\Z^+$ can be written as
$x^2+y^2+z^2+w^2\ (x,y,z,w\in\Z)$ with $ax+by+cz+dw\in\{8^k:\ k\in\N\}$, if and only if $(a,b,c,d)$ is among the quadruples
$$\gather(7,3,2,1),\ (7,5,2,1),\ (8,4,2,1),\ (8,4,3,2),\ (8,5,4,2),\ (8,6,2,1),
\\ (8,6,3,2),\ (8,6,5,1),\ (8,6,5,2),\ (9,8,7,4),\ (10,4,3,1),\ (12,4,3,1).
\endgather$$

{\rm (ii)} Let $a\in\Z^+$ and let $b,c,d\in\N$ with $a\gs b\gs c\gs d$ and $2\nmid \gcd(a,b,c,d)$. Then
any $n\in\Z^+$ can be written as
$x^2+y^2+z^2+w^2\ (x,y,z,w\in\Z)$ with $ax+by+cz+dw\in\{2\times 8^k:\ k\in\N\}$, if and only if $(a,b,c,d)$ is among the quadruples
$$\gather(7,3,2,1),\ (7,5,2,1),\ (8,4,2,1),\ (8,4,3,2),\ (8,5,4,2),\ (8,6,2,1),
\\ (8,6,3,2),\ (8,6,5,2),\ (9,5,4,2),\ (14,3,2,1),\ (14,5,2,1),\ (14,7,3,2).
\endgather$$
\endproclaim

\proclaim{Conjecture 4.7} {\rm (i) (24-Conjecture)} Any $n\in\N$ can be written as $x^2+y^2+z^2+w^2$ with $x,y,z,w\in\N$
such that both $x$ and $x+24y$ are squares.

{\rm (ii)} Each $n\in\N$ can be written as $x^2+y^2+z^2+w^2$ with $x,y,z,w\in\N$
such that both $x$ and $49x+48(y-z)$ are squares. Also, any $n\in\N$ can be written as $x^2+y^2+z^2+w^2$ with $x,y,z,w\in\N$
such that both $x$ and $121x+48(y-z)$ are squares.

{\rm (iii)} Every $n\in\N$ can be written as $x^2+y^2+z^2+w^2$ with $x,y,z,w\in\N$
such that both $x$ and $-7x-8y+8z+16w$ are squares.

{\rm (iv)} Each $n\in\N$ can be written as $x^2+y^2+z^2+w^2$ with $x,y,z,w\in\N$ and $x\eq y\pmod2$
such that both $x$ and $x^2+62xy+y^2$ are squares.
\endproclaim
\Remark\ 4.4. See [S, A281976, A281977, A281980, A282013, A282014, A282226 and A282463] for related data. We verified the 24-Conjecture for all $n=1,\ldots,10^7$, and Qing-Hu Hou extended the verification for $n$ up to $10^{10}$.
The author would like to offer 2400 US dollars as the prize for the first solution of the 24-Conjecture.
We verified parts (ii) and (iii) of Conjecture 4.7 for $n$ up to $10^7$ and $10^6$ respectively, and later Qing-Hu Hou extended
the verification of parts (ii) and (iii) of Conjecture 4.7 for $n$ up to $10^9$ and $10^8$ respectively.
Note that any $n\in\N$ is the sum of a fourth power and three squares as proved by the author [S17b].

\proclaim{Conjecture 4.8}
{\rm (i)} Any $n\in\Z^+$ can be written as $x^4+y^2+z^2+w^2$ with $x,y,z\in\N$ and $w\in\Z^+$
such that $9y^2-8yz+8z^2$ is a square. Also, any $n\in\N$ can be written as $4x^4+y^2+z^2+w^2$ with $x,y,z,w\in\N$
such that $79y^2-220yz+205z^2$ is a square.

{\rm (ii)} Each $n\in\Z^+$ can be written as $x^2+y^2+z^2+w^2$ with $x,y,z\in\Z$ and $w\in\Z^+$ such that both $2x+y$ and $2x+z$ are squares.

{\rm (iii)} Any $n\in\N$ can be written as $x^2+y^2+z^2+w^2$ with $x,w\in\N$ and $y,z\in\Z$ such that both $x+2y$ and $z+2w$ are squares.
Also, any $n\in\N$ can be written as $x^2+y^2+z^2+w^2$ with $x,y,z\in\Z$ and $w\in\N$ such that both $x+3y$ and $z+3w$ are squares.
\endproclaim
\Remark\ 4.5. See [S, A282933, A282972, A283170, A283196, A283204 and A283205] for related data.
The author [S17b] proved that any $n\in\N$ can be written as $x^2+y^2+z^2+w^2$ with $x,y,z,w\in\Z$ such that
$x+2y$ is a square, and that any $n\in\N$ can be written as $4x^4+y^2+z^2+w^2$ with $x,y,z,w\in\Z$.

\proclaim{Conjecture 4.9} {\rm (i)} Any $n\in\N$ can be written as $x^2+y^2+z^2+w^2$ with $x,y,z,w\in\N$
such that $x$ or $y$ is a square, and  $x-y$ is also a square.

{\rm (ii)} Any positive integer can be written as $x^2+y^2+z^2+w^2$ with $x,y,z,w\in\N$
such that $x+3y+5z$ is a positive square, and one of $2x,y,z$ (or $3x,y,z$) is a square.

{\rm (iii)} Any $n\in\Z^+$ can be written as $x^2+y^2+z^2+w^2$ with $x,y,z\in\N$ and $w\in\Z^+$
such that $(3x)^2+(4y)^2+(12z)^2$ is a square, and also one of $z,2z,3z$ is a square.
Also, each $n\in\Z^+$ can be written as $x^2+y^2+z^2+w^2$ with $x,y,z\in\N$ and $w\in\Z^+$
such that $(12x)^2+(15y)^2+(20z)^2$ is a square, and also one of $x,y,z$ is a square,
and any $n\in\N$ can be written as $x^2+y^2+z^2+w^2$ with $x,y,z,w\in\N$
such that $(12x)^2+(21y)^2+(28z)^2$ is a square, and also one of $x,2y,z$ is a square,
\endproclaim
\Remark\ 4.6. See [S, A281975, A300708, A300139, A300666, A300667, A300712, A300751, A300752, A300791, A300844, A300908] for related data or similar conjectures.

\proclaim{Conjecture 4.10} Any $n\in\N\sm\{71,85\}$ can be written as $x^2+y^2+z^2+w^2$ with $x,y,z,w\in\Z$
such that $9x^2+16y^2+24z^2+48w^2$ is a square.
\endproclaim
\Remark\ 4.7. We have verified this for $n$ up to $2\times10^5$. See [S, A281659] for related data and other similar conjectures.

\proclaim{Conjecture 4.11} {\rm (i)} Let $a,b\in\Z^+$ with $\gcd(a,b)$ squarefree. Then, each $n\in\N$
can be written as $x^2+y^2+z^2+w^2$ with $x,y,z,w\in\N$ such that $ax^3+b(y-z)^3$ is a square,
if and only if $(a,b)$ is among the ordered pairs
$$\gather (1,1),\ (1,9),\ (2,18),\ (8,1),\ (9,5),\ (9,8),
\\ (9,40),\ (16,2),\ (18,16),\ (25,16),\ (72,1).
\endgather$$

{\rm (ii)} Let $a,b\in\Z^+$ with $a\ls b$ and $\gcd(a,b)$ squarefree. Then, each $n\in\N$
can be written as $x^2+y^2+z^2+w^2$ $(x,y,z,w\in\Z)$ with $ax^3+by^3$ a square,
if and only if $(a,b)$ is among the ordered pairs
$$\gather (1,2),\ (1,8),\ (2,16),\ (4,23),\ (4,31),
\\ (5,9),\ (8,9),\ (8,225),\ (9,47),\ (25,88),\ (50,54).
\endgather$$
\endproclaim
\Remark\ 4.8. See [S, A282863 and A283617] for related data or similar conjectures.

\proclaim{Conjecture 4.12} {\rm (i)} Any $n\in\Z^+$ can be written as $x^2+y^2+z^2+3w^2$ $(x,y,z,w\in\Z)$ with $P(x,y,z,w)=1$,
whenever $P(x,y,z,w)$ is among the polynomials
$$\gather y+z+w,\ 2y+z+w,\ 2y+z+3w,
\ 2x+2y+z+w,
\\ 2x+2y+z+3w,\ 4x+2y+z+dw\ (d=1,3,5).\endgather$$

{\rm (ii)} Any $n\in\Z^+$ can be written as $x^2+y^2+2z^2+3w^2$ $(x,y,z,w\in\Z)$ with $P(x,y,z,w)=1$, whenever
$P(x,y,z,w)$ is among the polynomials
$$x+2y+w,\ y+z+w,\ y+2z+w,\ y+2z+3w,\ x+2y+2z+w,\ x+2y+2z+3w.$$

{\rm (iii)} Any $n\in\Z^+$ can be written as $x^2+y^2+3z^2+4w^2$ $(x,y,z,w\in\Z)$ with $y+z+2w=1$.
Also,  any $n\in\Z^+$ can be written as $x^2+y^2+2z^2+5w^2$ $(x,y,z,w\in\Z)$ with $y+2z+w=1$.
\endproclaim

\proclaim{Conjecture 4.13} {\rm (i)} Any $n\in\Z^+$ can be written as $x^2+y^2+z^2+2w^2$ $(x,y,z,w\in\Z)$ with $P_1(x,y,z,w)=1$,
whenever $P_1(x,y,z,w)$ is among the polynomials
$$\gather x+2y+3z,\ x+2y+5z,\ x+3y+4z,\ y+3z+2w,\ y+3z+4w,\ 2y+z+w,
\\x+y+2z+2w,\ x+2y+2z+2w,\ x+2y+3z+dw\ (d=1,2,4),
\\ x+2y+5z+2w,\ x+2y+5z+6w,\ x+3y+4z+2w,\ x+3y+4z+4w.
\endgather$$

{\rm (ii)} Any $n\in\Z^+$ can be written as $x^2+y^2+z^2+2w^2$ $(x,y,z,w\in\Z)$ with $P_2(x,y,z,w)=2$,
whenever $P_2(x,y,z,w)$ is among the polynomials
$$\gather x+y+2z+2w,\ x+y+2z+6w,\ x+2y+3z+2w,
\\  x+2y+3z+6w,\ x+2y+4z+4w,\ x+2y+5z+2w,\ 3x+3y+2z+2w.
\endgather$$

{\rm (iii)} Any $n\in\Z^+$ can be written as $x^2+y^2+z^2+2w^2$ $(x,y,z,w\in\Z)$ with $P_3(x,y,z,w)=3$,
whenever $P_3(x,y,z,w)$ is among the polynomials
$$x+2y+3z+2w,\ x+2y+3z+4w,\ x+2y+3z+6w.$$
\endproclaim

\proclaim{Conjecture 4.14} Any $n\in\Z^+$ can be written as $x^2+y^2+z^2+2w^2$ with $x,y,z,w\in\N$ and $x+2y\in\{4^k:\ k\in\N\}$.
Also, we may replace $x+2y$ by $y-z+3w$ (or $y+2z-w$).
\endproclaim
\Remark\ 4.9. It is easy to show that any $n\in\Z^+$ can be written as $x^2+y^2+z^2+2w^2$ with $x,y,z,w\in\Z$ and $x+2y=1$.
\medskip

Now we pose a conjecture similar to the 1-3-5 conjecture of Sun [S17b, Conjecture 4.3(i)].

\proclaim{Conjecture 4.15 {\rm (1-2-3 Conjecture)}} {\rm (i)} Any $n\in\N$ can be written as $x^2+y^2+z^2+2w^2$ with $x,y,z,w\in\N$ such that $x+2y+3z$ is a square.

{\rm (ii)} For each $n\in\Z^+$ we can write $n^2$ as $x^2+y^2+z^2+w^2$ with $x,y,z,w\in\N$ such that $x+2y+3z\in\{4^k:\ k\in\N\}$.
\endproclaim

\Remark\ 4.10. See [S, A275344 and A299924] for related data. Each of the numbers
 $$0,\ 1,\ 3,\ 5,\ 7,\ 14,\ 15,\ 16,\ 25,\ 30,\ 84,\ 169,\ 225$$
 has a unique representation $x^2+y^2+z^2+2w^2\ (x,y,z,w\in\N)$ with $x+2y+3z$ a square. For example,
$$\align 33 = &1^2 + 0^2 + 0^2 + 2\times4^2\ \ \t{with}\ 1 + 2\times0 + 3\times0 = 1^2,
\\84 = &4^2 + 6^2 + 0^2 + 2\times4^2\ \ \t{with}\ 4 + 2\times6 + 3\times0 = 4^2,
\\169 = &10^2 + 6^2 + 1^2 + 2\times4^2\ \ \t{with}\ 10 + 2\times6 + 3\times1 = 5^2,
\\225 = &10^2 + 6^2 + 9^2 + 2\times2^2\ \ \t{with}\ 10 + 2\times6 + 3\times9 = 7^2.
\endalign$$
Also, for each $n\in\{1,2,3,7,11,13,14,17,49,61\}$ there is a unique way to write $n^2$ as $x^2+y^2+z^2+w^2$
with $x,y,z,w\in\N$ and $x+2y+3z\in\{4^k:\ k\in\N\}$. For example,
$$\align 11^2=&2^2+1^2+4^2+10^2\ \ \t{with}\ 2+2\times1+4\times3=4^2,
\\49^2=&22^2+3^2+12^2+42^2\ \ \t{with}\ 22+2\times3+3\times12=4^3.
\endalign$$
We conjecture that if $a,b,c$ are positive integers with $\gcd(a,b,c)$ squarefree and any positive square can be written as
$x^2+y^2+z^2+w^2$ with $x,y,z,w\in\N$ and $ax+by+cz\in\{4^k:\ k\in\N\}$ then we must have $\{a,b,c\}=\{1,2,3\}$. \medskip

\proclaim{Conjecture 4.16} Let $\da\in\{0,1\}$.

{\rm (i)} For any integer $n>\da$, we can
write $n^2$ as $x^2+y^2+z^2+w^2$
with $x,y,z,w\in\N$ and $\{x,4x-3y\}\se\{2^{2k+\da}:\ k\in\N\}$.

{\rm (ii)} For any $n\in\Z^+$ we can write $2n^2=x^2+y^2+z^2+w^2$ with $x,y,z,w\in\N$ such that
$x+3y+5z+15w\in\{2^{2k+\da}:\ k\in\N\}$.
\endproclaim
\Remark\ 4.11. We have verified part (i) for $n$ up to $10^7$. See [S, A300219, A299537, A299794, A300360, A300396, A301891] for related data or similar conjectures.
\medskip

\Ack. The author would like to thank the referee for his helpful comments on Lemma 3.4 and its proof.

\enddocument